\documentclass[[a4paper,11pt]{article}
\usepackage{amsmath}
\usepackage{amsfonts}
\usepackage{amssymb}
\usepackage{epsfig}
\usepackage{latexsym}
\usepackage{graphicx}
\usepackage{epstopdf}
\usepackage{float}
\usepackage{mathrsfs}
\usepackage{amsmath}
\usepackage{setspace}
\usepackage{caption}
\usepackage{multirow}

\newtheorem{theorem}{\bf Theorem}[section]
\newtheorem{lemma}[theorem]{\bf Lemma}

\newtheorem{coro}[theorem]{\bf Corollary}
\newtheorem{example}[theorem]{\bf Example}

\newenvironment{remark}{\noindent{\bf Remark}}{\vspace{2ex}}

\newcommand{\Span} {{\rm span }}

\def\mA{\mathbb{A}^2}
\def\D{\bf D}
\def\qed{\quad \vrule height7.5pt width4.17pt depth0pt}

\begin{document}

\title{Rational Approximation in the Bergman Spaces}
\author{Wei QU\thanks{Department of Information Technology, 
Macau University of Science and Technology. Email: 376928001@qq.com.\newline \indent { Department of Information Technology, 
Macau University of Science and Technology. Email: pdang@must.edu.mo.}}  ,  Pei DANG}
\maketitle

\begin{abstract} It is known that adaptive Fourier decomposition (AFD) offers efficient rational approximations to functions in the classical Hardy $H^2$ spaces with significant applications. This study aims at rational approximation in Bergman, and more widely, in weighted Bergman spaces, the functions of which have more singularity than those in the Hardy spaces. Due to lack of an effective inner function theory, direct adaptation of the Hardy-space AFD is not performable. We, however, show that a pre-orthogonal method, being equivalent to AFD in the classical cases, is available for all weighted Bergman spaces. The theory in the Bergman spaces has equal force as AFD in the Hardy spaces. The methodology of approximation is via constructing the rational orthogonal systems of the Bergman type spaces, called Bergman space rational orthogonal (BRO) system, that have the same role as the Takennaka-Malmquist (TM) system in the Hardy spaces. Subsequently, we prove a certain type direct sum decomposition of the Bergman spaces that reveals the orthogonal complement relation between the span of the BRO system and the zero-based invariant spaces. We provide a sequence of examples with different and explicit singularities at the boundary along with a study on the inclusion relations of the weighted Bergman spaces. We finally present illustrative examples for effectiveness of the approximation.  \end{abstract}


\vspace *{0.5mm}

\section {Introduction}

The aim of this work is to propose an adaptive approximation method in the weighted Bergman $\mA_\alpha$ spaces analogous with the one established in the Hardy $H^2$ spaces of the unit disc and the upper-half complex plane (classical contexts) called adaptive Fourier decomposition (AFD) \cite{QWa}. The theory established in the Hardy spaces has direct applications in system identification and in signal analysis \cite{Mi1, Mi2, QianBook1}. It is well known that a function $f\in L^2(\partial \D)$ with Fourier expansion $f(e^{it})=\sum_{-\infty}^\infty c_ne^{int}$ has its Hardy space decomposition $f=f^++f^-,$ where
$f^+(e^{it})=\sum_{0}^\infty c_ne^{int}\in H^2_+(\partial \D),$ and $f^-(e^{it})=\sum_{-\infty}^{-1} c_ne^{int}\in H^2_-(\partial \D).$ This decomposition corresponds to the space direct sum decomposition
\begin{eqnarray}\label{DC} L^2(\partial \D)=H^2_+(\partial \D)\oplus H^2_-(\partial \D),\end{eqnarray}
where $H^2_\pm(\partial \D)$ are, correspondingly, the non-tangential boundary limits of the functions in the respective Hardy spaces inside and outside the unit disc, namely, $H^2_\pm (\D).$ The projection operators from $f$ to $H^2_\pm (\D)$ are, respectively, denoted as $P_\pm (f)=(1/2)(f\pm iHf\pm c_0).$ If $f$ is real-valued, due to the property $c_{-n}=\overline{c}_n,$ one has $f=2{\rm Re}f^+-c_0.$ The inverse operator of $P_+$ mapping $f^+$ to $f$ is invertible and bounded. Those suggest that analysis of the $L^2$ functions may be reduced to analysis of the functions in the corresponding Hardy spaces \cite{QianBook1}. Let $\{B_k\}$ be the rational orthogonal system defined by a sequence of numbers ${\bf a}=(a_1,\cdots,a_n,\cdots ),$ allowing repetition, in the unit disc of the form
\[ B_k(z)=\frac{\sqrt{1-|a_k|^2}}{1-\overline{a}_kz}\prod_{l=1}^{k-1}\frac{z-a_l}{1-\overline{a}_lz}.
\]
$\{B_k\}_{k=1}^\infty$ is an orthonormal system, also called Takenaka-Mulmquist or TM system, in $H^2_+(\D).$ The system may or may not be a basis of $H^2_+(\D)$ depending whether $\sum_{k=1}^\infty (1-|a_k|)=\infty$ or not. In rational approximation in Hardy spaces one cannot avoid the above defined rational orthonormal systems for they are results of the Gram-Schmidt orthogonalization applied to partial fractions with different poles. Related to TM systems is the Beurling-Lax type direct sum decomposition
\begin{eqnarray}\label{OR}
H^2(\D)=\overline{{\rm span}}\{B_k\}_{k=1}^\infty \oplus \phi H^2(\D),
\end{eqnarray}
where $\phi$ is the Blaschke product defined by the parameters $a_1,\cdots,a_k,\cdots$. The validity of the last relation rests on the condition $\sum_{k=1}^\infty (1-|a_k|)<\infty$ under which $\phi$ is well defined and the TM system is not a basis. Selecting parameters $a_k$'s under the maximal selection principle, regardless whether the corresponding TM system is a basis or not, the resulted AFD expansion
\[ f^+(z)=\sum_{k=1}^\infty \langle f,B_k\rangle B_k(z)\]
converges at fast pace \cite{QWa, QWang}. One extra property, but important in signal analysis, of the decomposition is that under the selection $a_1=0$ all $B_k$'s are of positive frequencies (boundary phase derivative). AFD thus gives rise to intrinsic positive frequency decomposition of $f^+$ and thus that of $f$ as well. Other related work can be seen in \cite{Baratchart,Wang, Mai}. In the analytic adaptive and positive-frequency approximation aspects the studies of Coifman et al. and Qian et al. merged together \cite{Coi1, Coi2}.

To extend the theory of the Hardy spaces to the Bergman spaces a decomposition of $L^2(\D)$ similar to (\ref{DC}) with bounded invertible projections does not seem to be possible. Simple computation shows that the projection from $L^2(\D)$ to $\mA(\D)$ is not one to one. The main task of this work is to explore the availability of (\ref{OR}) in the Bergman spaces case.
The remarkable difference between the Hardy space and, for instance, the standard Bergman space is that
\begin{eqnarray*}
\mA(\D)\ne\overline{\rm span}\{B_k\}_{k=1}^\infty \oplus {\it H}_{\bf a} \mA(\D),
\end{eqnarray*}
where $\{B_k\}_{k=1}^\infty$ is the Gram-Schmidt orthonormalization of a sequence of generalized reproducing kernels (see \S 2) corresponding to ${\bf a}$ that defines a Horowitz product $H_{\bf a}$ \cite{Duren} (see \S 2 and \S 5). On the other hand, we show that the direct sum decomposition holds with the Beurling type shift invariant subspace being replaced by the zero-based invariant subspace, that is
\begin{eqnarray*}
\mA(\D)=\overline{\rm span}\{B_k\}_{k=1}^\infty \oplus I_{\bf a}.
\end{eqnarray*}

\section {Preliminaries}

Denote by ${\bf D}$ the open unit disc in the complex plane ${\bf C}.$ We will first deal with the square integrable Bergman space of the open unit disc, $\mathbb{A}^2({\bf D}),$ namely,
$$ \mathbb{A}^2({\bf D})=\{f: {\D}\to {\bf C}\ | \ f\ {\rm is \ holomorphic \ in \ \D, \ and\ }\ \| f\|_{\mA({\bf D})}^2=\int_{\bf D} |f(z)|^2dA<\infty\},$$
where $dA$ is the normalized area measure on the unit disc: $dA=\frac{dxdy}{\pi}, z=x+iy.$ The space $\mA({\bf D})$ under the norm  $\| \cdot \|_{{\mA}({\bf D})}$ forms a Hilbert space with the inner product
$$ \langle f,g\rangle_{\mA({\bf D})}=\int_{\bf D} f(z)\overline{g(z)}dA.$$
In the sequel we often suppress the subscripts $\mA({\bf D})$ in the notations $\| \cdot \|_{\mA({\bf D})}$ and $\langle \cdot,\cdot\rangle_{\mA({\bf D})}$, and write them simply as $\|\cdot \|$ and $\langle \cdot,\cdot\rangle.$ We also write $\mA({\bf D})$ as $\mA$ \cite{Zhu,Duren}.

It is known that $\mA$ is a reproducing kernel Hilbert space with the reproducing kernel
$$ k_a(z)=\frac{1}{(1-\overline{a}z)^2}.$$
From the reproducing kernel property we have
$$ \| k_a\|^2= k_a(a)=\frac{1}{(1-|a|^2)^2}.$$
Thus, the normalized reproducing kernel, denoted as $e_a(z)$, is
$$ e_a(z)=\frac{k_a(z)}{\sqrt{k_a(a)}}=\frac{1-|a|^2}{(1-\overline{a}z)^2}.$$
We hence have, for any $f\in \mA,$
$$ \langle f,e_a\rangle = (1-|a|^2)f(a).$$

Let ${\bf a}$ be an infinite sequence with its elements in $\D,$ that is ${\bf a}=(a_1,\cdots,a_n, \cdots), a_n\in {\D}, n=1,2,\cdots .$ Here, and in the sequel, $a_n$'s are allowed to repeat. Denote by $l(a_n)$ as the number of repeating time of the $n$-th element $a_n$ of ${\bf a}_n$ in the $n$-tuple ${\bf a}_n=(a_1,\cdots,a_n).$
Denote by
\begin{eqnarray}\label{derivative}
\tilde{k}_{a_n}(z)=\left(\frac{d}{d\overline{w}}\right)^{(l(a_n)-1)}\left(k_w(z)\right)|_{w=a_n},
\end{eqnarray}
that is, the $(l(a_n)-1)$-th derivative of $k_{w}$ with respect to the variable $\overline{w},$ then replace $w$ by $a_n.$ We will call the sequence $\{\tilde{k}_{a_n}\}$ the {\it generalized reproducing kernels corresponding to}  $(a_1,\cdots,a_n)$.

We will be using the relation, for $f\in \mA,$
\begin{eqnarray}\label{beautiful} \langle f, \tilde{k}_{a_n}\rangle =f^{l(a_n)-1}(a_n). \end{eqnarray}
To show (\ref{beautiful}), taking the $(l-1)$ differentiation to the both sides of
the identity
\[ f(w)=\int_{\D} \frac{f(z)}{(1-w\overline{z})^2}dA,\]
we have
\[ f^{(l-1)}(w)=l!\int_{\D} \frac{\overline{z}^{l-1}f(z)}{(1-w\overline{z})^{l+1}}dA.\]
The exchange of the order of differentiation and integration is justified by the Lebesgue Dominated Convergence Theorem based on the fact that the integrand does not have singularity. Taking into account
\[ \overline{\left(\frac{d}{d\overline{w}}\right)^{l-1}}k_w(z)=\left(\frac{d}{dw}\right)^{l-1}\overline{k}_w(z),\]
we obtain (\ref{beautiful}).

In various contexts, as well as in general reproducing kernel Hilbert spaces, including the Bergman spaces, this relation was noted by a number of authors \cite{QWang, Q2}. As revealed in the following sections, various degrees of differentiations of the kernel function correspond to times of repetitions of the parameters $a_1,...,a_n,...$ Some literature for the clarity purpose deliberately avoid the repetition cases. To get the best approximations, repetition, however, is unavoidable and should be carefully treated.

We will be using the following functions, called Horowitz products, for $a\in \D$ and ${\bf a}_n=(a_1,\cdots,a_n)\in \D^n,$
\[ H_{\{a\}}(z)=\frac{|a|}{a}
\frac{a-z}{1-\overline{a}z}\left(2-\frac{|a|}{a}\frac{a-z}{1-\overline{a}z}\right)
\]
and
\[ H_{{\bf a}_n}= H_{\{a_1\}}\cdots H_{\{a_n\}}.\]
For ${\bf a}=(a_1,\cdots,a_n,\cdots),$
the condition
\begin{eqnarray}\label{H-product} \sum_{k=1}^\infty (1-|a_k|)^2<\infty \end{eqnarray}
is necessary and sufficient for convergence of the corresponding infinite Horowitz product
\[ H_{\bf a}(z)=\prod_{k=1}^\infty H_{\{a_k\}}(z),\qquad |z|<1.\]
Under the condition the convergence is
uniform for $|z|\leq r$ for any fixed $0<r<1$. It is noted that when
the infinite product converges it does not have to converge to a function in the Bergman space \cite{Duren}.

We will study zero-based invariant subspaces. Let $I$ be a closed subspace of $\mA.$  If further $zI\subset I,$ then we say that $I$ is an invariant subspace.

In the following we will concern the cases where $I$ corresponds to the set of functions that all vanish on an ordered set of points $A,$ where
the points of $A$ can repeat. Obviously such set $I$ is an invariant subspace, denoted as $I=I_A,$ and is called the $A$-zero based invariant subspace of $\mA:$
\[ I_A=\{f\in \mA\ :\ f \ {\rm has \ all \ points  \ in \ } A \ {\rm as\ its\ zeros}\}.\]
In the above definition when $a\in A$ repeats $l$ times, by saying that $f$ takes value zero at all $a\in A$ we mean that $f(a)=f'(a)=\cdots =f^{(l-1)}(a)=0.$ The fact $f\in I_A$ is also denoted $f(A)=0.$
We will be interested in the case $A=A_n={\bf a}_n=(a_1,\cdots , a_n)$
and $A={\bf a}=(a_1,\cdots, a_n, \cdots ).$ Below we will denote $\phi_a(z)=z-a.$  In the case
$I=I_{{\bf a}_1}=I_{\{a_1\}},$ we can show that $I_{\{a_1\}}=\phi_{a_1}\mA.$ First, it is easy to see
 $\phi_{a_1}\mA \subset I_{\{a_1\}}.$ On the other hand, if $f\in \mA$ and $f(a_1)=0,$ then $f(z)=\phi_{a_1}(z)g(z).$
It is easy to show $ g\in \mA.$ In fact, $g$ is analytic inside the disc $\D$, and hence bounded in $|z|\leq |a_1|+\epsilon <1.$
Then
\[ \int_{|z|\leq |a_1|+\epsilon} |g(z)|^2 dA(z)< \infty.\]
On the other hand,
\begin{eqnarray*}
 \int_{|a_1|+\epsilon<|z|<1} |g(z)|^2dA(z) &=&  \int_{|a_1|+\epsilon<|z|<1} \frac{|f(z)|^2}{|z-a_1|^2}dA(z)\\
 &\leq & \frac{1}{\epsilon^2}  \int_{\D} |{f(z)}|^2dA(z)\\
 &<& \infty.\end{eqnarray*}
Therefore, $g\in \mA,$ and $I_{{\bf a}_1}\subset \phi_{a_1}\mA.$
Inductively, we can show $I_{{\bf a}_n}=\phi_{a_1}\cdots \phi_{a_n}\mA.$
It is easy to see that $H_{{\bf a}_n}\mA$ is identical with $\phi_{a_1}\cdots \phi_{a_n}\mA.$ Since the zero-based
invariant spaces are closed subspaces, they themselves are reproducing kernel Hilbert spaces.
When ${\bf a}$ is the zero set of some function $f\in \mA,$ then there holds the relation (\ref{H-product}).
An infinite sequence ${\bf a}$ satisfying condition (\ref{H-product}), however,
 does not imply that it is the zero set of a function in the Bergman space.
 Necessary and sufficient conditions for ${\bf a}$ being the zero set of some
 Bergman space function is so far an open problem \cite{C,prob,D}.

\section {Rational Orthogonal System of $\mA$}

In this section we study rational orthogonal systems, or Takenaka-Malmquist systems as an alternative terminology,  in the Bergman space.
The following theorem is crucial. The scope of the following theorem is as the same as that in \cite{Zhupaper}, except that multiples of the parameters are allowed. This generalization turns to be important for achieving the best approximation with the pre-orthogonal algorithm given in the next section (also see \cite{Q2}).

\begin{theorem}\label{th1} Let $(a_1,\cdots, a_n, \cdots)$ be an infinite sequence in $\D,$ where each $a_n$ in the sequence is allowed to repeat. Let the orthonormalization of $(\tilde{k}_{a_1},\cdots,\tilde{k}_{a_n})$ be denoted
$(B_1,...,B_n).$ Then for any $n$, (i) the reproducing kernel of the zero-based invariant subspace $I_{A_n}$ is
$$K_{A_n}(z,w)=\sqrt{\|k_w\|^2-\sum_{k=1}^n\langle k_w, B_k\rangle|^2}B_{n+1}^w,\ K_{A_0}(z,w)=K(z,w),$$
where $(B_1,...,B_n,B_{n+1}^w)$ is the orthonormalization of $(\tilde{k}_{a_1},\cdots,\tilde{k}_{a_n},k_w), w$ is none of $a_1,\cdots,a_n;$ (ii) if $a_{n+1}$ coincides with some $a_k, k=1,\cdots,n,$ that is, $l(a_{n+1})>1,$ then $B_{n+1}^{a_{n+1}}$ is defined through
\begin{eqnarray*}
\lim_{w\to a_{n+1}}B_{n+1}^w&=& \frac{\tilde{k}_{a_{n+1}}(z)-\sum_{k=1}^n \langle \tilde{k}_{a_{n+1}}, B_k\rangle B_k (z)}{\|\tilde{k}_{a_{n+1}}-\sum_{k=1}^n \langle \tilde{k}_{a_{n+1}}, B_k\rangle B_k \|}\\
&=&\frac{\tilde{k}_{a_{n+1}}(z)-\sum_{k=1}^n \langle \tilde{k}_{a_{n+1}}, B_k\rangle B_k (z)}{\sqrt{\|k_w\|^2-\sum_{k=1}^n|\langle \tilde{k}_{a_{n+1}}, B_k\rangle|^2}},
\end{eqnarray*}
and thus $B_{n+1}\triangleq B_{n+1}^{a_{n+1}}.$
\end{theorem}

\noindent{\bf Proof.}
(i) Let $w$ be none of $a_1,\cdots,a_n.$ Under the Gram-Schmidt (G-S) orthonormalization process, $B_{n+1}^w$ is defined as
\begin{eqnarray}\label{n=n} B_{n+1}^w(z)=\frac{k_{w}(z)-\sum_{k=1}^n \langle k_{w}, B_k\rangle B_k (z)}{\|k_{w}-\sum_{k=1}^n \langle k_w, B_k\rangle B_k \|}=\frac{k_{w}(z)-\sum_{k=1}^n \langle k_{w}, B_k\rangle B_k (z)}{\sqrt{\|k_w\|^2-\sum_{k=1}^n|\langle k_w, B_k\rangle|^2}}.
\end{eqnarray}
We are to show that
\[ K_{A_n}(\cdot,w)=\sqrt{\|k_w\|^2-\sum_{k=1}^n|\langle k_w, B_k\rangle|^2}B_{n+1}^w(z)=k_{w}-\sum_{k=1}^n \langle k_{w}, B_k\rangle B_k (z)\]
is the reproducing kernel of $I_{A_n}.$ To this aim we first show (1) $K_{A_n}(\cdot,w)$ itself belongs to the space $I_{A_n};$ and then show (2) For $f=\phi_{a_1}\cdots\phi_{a_n}g, g\in \mA,$ there holds $\langle f, K_{A_n}(\cdot,w)\rangle =f(w).$

Now we show (1). From the construction of $K_{A_n}(\cdot,w)$ we know that it is orthogonal with all $B_1,\cdots B_n,$ and hence orthogonal with all $\tilde{k}_{a_1},\cdots,\tilde{k}_{a_n}.$ Let $a$ be a complex number appearing in the sequence $(a_1,\cdots,a_n)$ with the repetition number $l.$ Then $k_a, \frac{\partial}{\partial \overline{w}}k_a,\cdots,\left(\frac{\partial}{\partial \overline{w}}\right)^{l-1}k_a$ will appear in  $\tilde{k}_{a_1},\cdots,\tilde{k}_{a_n}.$ By invoking (\ref{beautiful}), the function $K_{A_n}(\cdot,w)$ being orthogonal with $k_a, \frac{\partial}{\partial \overline{w}}k_a,\cdots,\left(\frac{\partial}{\partial \overline{w}}\right)^{l-1}k_a$ means that $K_{A_n}(a,w)=\frac{\partial}{\partial z}K_{A_n}(a,w)=\cdots =\left(\frac{\partial}{\partial z}\right)^{l-1}K_{A_n}(a,w)=0.$ The latter assertion further implies that $K_{A_n}(z,w)$ has the multiplicative factor $(z-a)^l$ and thus has the factor $\phi_a^l.$ By applying this argument for all $a$ appearing in the sequence $(a_1,\cdots,a_n)$ we obtain that for each $w$, $K_{A_n}(\cdot,w)$ is in $I_{A_n}.$

Next we show (2), that is, for $f=\phi_{a_1}\cdots\phi_{a_n}g, g\in \mA,$ there holds $\langle f, K_{A_n}(\cdot,w)\rangle =f(w).$ We first claim that
\begin{eqnarray}\label{claim}
\langle \phi_{a_1}\cdots\phi_{a_n}g, \tilde{k}_{a_k}\rangle=0,\ \ \ \ \ k=1,2,\cdots,n.
\end{eqnarray}
 Let $a_k$ repeat $l$ times in $(a_1,\cdots,a_n)\ $ and $l(a_k)$ times in $(a_1,\cdots,a_k).$ We always have $l\geq l(a_k).$ In such notation the factor$(\phi_{a_k})^l$ lies in $\phi_{a_1}\cdots\phi_{a_n}$ and for some $h\in \mA,\ \phi^l_{a_k}h=\phi_{a_1}\cdots\phi_{a_n}g.$ By (\ref{beautiful}), we have
\begin{eqnarray*}
  \langle \phi_{a_1}\cdots\phi_{a_n}g, \tilde{k}_{a_k}\rangle &=& \langle \phi^l_{a_k} h, \left(\frac{\partial}{\partial \overline{w}}\right) ^{l(a_k)-1}k_{w} \rangle|_{w=a_k}\\
  &=& (\phi^l_{a_k}h)^{(l(a_k)-1)}(a_k)\\
  &=& 0.
\end{eqnarray*} Since each $B_k$ is a linear combination of $(\tilde{k}_{a_1},\cdots,\tilde{k}_{a_k}),$ we have, owing to (\ref{claim}),
\begin{eqnarray*}
 \langle\phi_{a_1}\cdots\phi_{a_n}g, K_n(\cdot,w)\rangle&=& \langle\phi_{a_1}\cdots\phi_{a_n}g, k_{w}-\sum_{k=1}^n \langle k_{w}, B_k\rangle B_k (z)\rangle\\
 &=&\langle \phi_{a_1}\cdots\phi_{a_n}g, k_w\rangle +\sum_{k=1}^n c_k \langle \phi_{a_1}\cdots\phi_{a_n}g, \tilde{k}_{a_k}\rangle\\
 &=& \phi_{a_1}(w)\cdots\phi_{a_n}(w)g(w)\\
 &=& f(w).
\end{eqnarray*}
Therefore, $K_n(z,w)$ is the reproducing kernel of $I_{A_n}.$
The proof of (2) is complete.

 Now we prove (ii) that treats $l=l(a_{n+1})>1.$ This means that, with a little abuse of notation, the $l-1$ terms $k_{a_{n+1}}, \frac{\partial}{\partial \overline{w}}k_{a_{n+1}},\cdots, $ $\left(\frac{\partial}{\partial \overline{w}}\right)^{(l-2)}k_{a_{n+1}}$ as functions of the $z$ variable have appeared in the sequence $(\tilde{k}_{a_1}, \cdots,\tilde{k}_{a_{n}}).$ Therefore, the function
\[ k_{a_{n+1}}+\frac{\frac{\partial}{\partial \overline{w}}k_{a_{n+1}}}{1!}(\overline{w}-\overline{a}_{n+1})+\cdots
+\frac{\left(\frac{\partial}{\partial \overline{w}}\right)^{(l-2)}k_{a_{n+1}}}{(l-2)!}(\overline{w}-\overline{a}_{n+1})^{l-2},\]
as the order-$(l-2)$ Taylor expansion of the function $k_w(z)$ with respect to the variable $\overline{w}$ at $\overline{w}=\overline{a}_{n+1},$ is already in the linear span of $B_1,\cdots, B_n.$ Denoting the Taylor expansion by $T_{l-2}(k_w, a_{n+1}),$ the last mentioned assertion amounts to the relation
 \begin{eqnarray}\label{insert1} T_{l-2}(k_w, a_{n+1})-\sum_{k=1}^n\langle T_{l-2}(k_w, a_{n+1}),B_k\rangle B_k=0.\end{eqnarray}
For $w$ being different from all $a_k, k=1,\cdots, n,$ we have, owing to (\ref{n=n}),
\begin{eqnarray}\label{into1} \frac{k_{w}(z)-\sum_{k=1}^n \langle k_w, B_k\rangle B_k (z)}{\|k_{w}-\sum_{k=1}^n \langle k_w, B_k\rangle B_k \|}=B_{n+1}^w(z),\end{eqnarray}
Inserting (\ref{insert1}) into (\ref{into1}), and dividing by $(\overline{w}-\overline{a}_{n+1})^{l-1}$ and
$|\overline{w}-\overline{a}_{n+1}|^{l-1}$ to the numerator and the denominator, respectively, we have
\begin{eqnarray}\label{inserted} e^{-(l-1)\theta}\frac{\frac{k_{w}(z)-T_{l-2}(k_{w}, a_{n+1})(z)}{(\overline{w}-\overline{a}_{n+1})^{l-1}}-\sum_{k=1}^n \langle \frac{k_{w}-T_{l-2}(k_{w}, a_{n+1})(z)}{\overline{w}-\overline{a}_{n+1})^{l-1}}, B_k\rangle B_k (z)}{\|\frac{k_{w}(z)-T_{l-2}(k_{w}, a_{n+1})(z)}{(\overline{w}-\overline{a}_{n+1})^{l-1}}-\sum_{k=1}^n \langle \frac{k_{w}(z)-T_{l-2}(k_{w}, a_{n+1})(z)}{(\overline{w}-\overline{a}_{n+1})^{l-1}}, B_k\rangle B_k \|}=B_{n+1}^w(z),\end{eqnarray}
where  $\overline{w}-\overline{a}_{n+1}=|\overline{w}-\overline{a}_{n+1}|e^{i\theta}.$
Letting $w\to a_{n+1}$, keeping the direction $e^{i\theta},$ and using the Lagrange type remainder, we obtain
\[ e^{-(l-1)\theta}\frac{\tilde{k}_{a_{n+1}}(z)-\sum_{k=1}^n \langle \tilde{k}_{a_{n+1}}, B_k\rangle B_k (z)}{\|\tilde{k}_{a_{n+1}}-\sum_{k=1}^n \langle \tilde{k}_{a_{n+1}}, B_k\rangle B_k \|}=B_{n+1}(z).\]
This shows that the G-S orthonormalization of $(\tilde{k}_{a_1},\cdots,\tilde{k}_{a_{n+1}})$ is
$(B_1,\cdots,B_{n+1}).$ The proof of (ii) is complete.
\rightline{\qed}

We call $(B_1,\cdots, B_n, \cdots )$ the rational orthonormal system or the Takenaka-Malmquist system (TM system) of the Bergman space associated with the infinite sequence $(a_,\cdots,a_n,\cdots )$ or, equivalently, the generalized reproducing kernels sequence $(\tilde{k}_{a_1},\cdots,\tilde{k}_{a_n},\cdots),$ where for each $n$, $(B_1,\cdots, B_n)$ is the Gram-Schmidt orthonormalization of  $(\tilde{k}_{a_1},\cdots,\tilde{k}_{a_n}).$
For fixed $a_1,\cdots,a_n,$ the correspondence between $a_{n+1}$ and $B_{n+1}$ is unique up to a multiplicative
constant of complex modular $1.$

\section {Pre-Orthogonal Adaptive Fourier Decomposition of the Bergman Space $\mA$}

The purpose of this section is to introduce a machinery that gives rise to fast rational approximation to functions in the Bergman space. The machinery is an adaptation of a general method called pre-orthogonal adaptive Fourier decomposition or POAFD for certain Hilbert spaces possessing the boundary vanish property \cite{Q2}. The major work of this section is to verify the boundary vanishing condition of the Bergman space. It is the pre-orthogonal process that requires repeating selections of the parameter $a,$ when necessary, and therefore the necessity of the derivatives of the reproducing kernel.

Let $f\in \mA.$ Denoting $f=f_1.$ For the normalized reproducing kernel $e_a$ for any $a\in \D,$ we have the identity
\[ f(z)=\langle f_1, e_{a_1}\rangle e_{a_1}(z) + f_2(z),\]
where $f_2=f_1-\langle f_1, e_{a_1}\rangle e_{a_1}$ is orthogonal with $e_{a_1}.$ The last assertion is proved through the Hilbert space projection property. We also have
\[ \|f\|^2=|\langle f_1,e_{a_1}\rangle|^2+\|f_2\|^2.\]
The strategy is to maximize the value of $|\langle f_1,e_{a_1}\rangle|^2$ through selections of $a_1,$ and thus to minimize the energy of $f_2.$ Despite of the fact that ${\bf D}$ is an open set, we have

\begin{lemma}\label{premaximal selection}
For any $f\in \mA$ there exists $a\in {\bf D}$ such that
\[ |\langle f,e_{a}\rangle|^2=\max \{ |\langle f,e_{b}\rangle|^2\ :\ b\in {\bf D}\}.\]
\end{lemma}

\noindent{\bf Proof}
The Cauchy-Schwarz inequality gives
\[|\langle f,e_b\rangle|\leq \|f\|.\]
Therefore $ |\langle f,e_b\rangle| $ has a finite upper bound.

For $\epsilon >0,$ since polynomials are dense in the Bergman space, there exists a polynomial $g$ such that
\[ \|f-g\|\leq \epsilon.\]
The Cauchy-Schwarz inequality gives
\[ |\langle f,e_b\rangle|\leq |\langle g,e_b\rangle|+\epsilon=(1-|b|^2)|g(b)|+\epsilon.\]
The last quantity tends to zero as $|b|\to 1.$ Therefore $|\langle f,e_b\rangle|$ attains the maximum value at an interior point.  \qed

For $f_1\in \mA$ choose $a_1$ such that
\[  |\langle f_1,e_{a_1}\rangle|^2=\max \{ |\langle f_1,e_{b}\rangle|^2\ :\ b\in {\bf D}\}.\]
Fixing such $a_1,$ we have, by denoting $e_{a_1}=B_1,$
 \begin{eqnarray}\label{f3} f(z)=\langle f_1, B_1\rangle B_1(z) + \langle f_2, B_2^b\rangle B_2^b(z) + f_3,\end{eqnarray}
 where $(B_1,B_2^b)$ is the G-S orthonormalization of $(B_1,\tilde{k}_{b}).$ In the end of the argument we will not exclude the case $b=a_1,$ but for the time being we assume $b\ne a_1.$ By invoking the formula (\ref{n=n}) and the orthogonality between $f_2$ and $B_1,$ we have

 \begin{eqnarray*} \langle f_2,B_2^b\rangle = \langle f_2, \frac{e_b}{\sqrt{1-|\langle e_b,B_1\rangle|^2}} \rangle \end{eqnarray*}

 The corresponding integral, first, has a finite upper bound independent of $b\in {\bf D}.$ We will show that the norm of the complex number depending on $b$ in the above equality
 reaches its maximum at an interior point $a_2\in {\bf D}.$ To this end it suffices to show that
 the quantity tends to zero as $|b|\to 1.$ The quantity can be decomposed into two parts, as
 \[ \langle f_2, \frac{e_b}{\sqrt{1-|\langle e_b,B_1\rangle|^2}} \rangle=\langle f_2, \frac{e_b}{\sqrt{1-|\langle e_b,B_1\rangle|^2}}- e_b \rangle + \langle f_2, e_b \rangle= I_1(b)+I_2(b).\]
 \def\bD{\bf D}
 We have
 \begin{eqnarray}\label{using} |I_1(b)|\leq \int_{\bD} |f_2(z)||e_b(z)|\left(\frac{1}{\sqrt{1-|\langle e_b,B_1\rangle|^2}}-1\right)dA(z).\end{eqnarray}
Owing to the reproducing kernel property, for any $\epsilon>0,$ if $|b|$ is close to $1-,$ then
 \[|\langle e_b,B_1\rangle|=\frac{(1-|b|^2)(1-|a_1|^2)}{|1-\overline{a}_1b|^2}\leq \epsilon.\]
 Using the last obtained estimate in (\ref{using}), and then invoking the Cauchy-Schwarz inequality, we obtain
 \[ \lim_{|b|\to 1-}I_1(b)=0.\]
 The maximal selection Lemma \ref{premaximal selection} implies $\lim_{|b|\to 1-}I_2(b)=0 .$ To summarize, we have
 \begin{eqnarray}\label{same reasoning} \lim_{|b|\to 1-} |\langle f_2,B_2^b\rangle|=0.\end{eqnarray}
 The process of reaching and attaining the maximal value of $|\langle f_2,B_2^b\rangle|$ can always be made from a sequence $\{b_k\}_{k=1}^\infty$ such that $\lim_{k\to \infty}b_k=a_2$ and for each $k, b_k\ne a_2.$ If $a_2\ne a_1,$ then one simply has $(B_1,B_2)=(B_1,B^{a_2}_2),$ the latter being the G-S orthonormalization of $(k_{a_1}, k_{a_2}).$ If, finally, $a_2=a_1,$ then according to Theorem \ref{th1}, one also has $(B_1,B_2)=(B_1,B^{a_2}_2),$ but the latter is the G-S orthonormalization of $(k_{a_1}, \frac{\partial}{\partial \overline{w}}k_{a_2}).$ By the definition of generalized reproducing kernels in relation to the ordered pairs $(a_1, a_2),$ both cases will be the G-S orthonormalization of $(\tilde{k}_{a_1}, \tilde{k}_{a_2}).$  Fixing such a maximal selection $a_2$ we have that $f_3$ in (\ref{f3}) is orthogonal to $B_1$ and $B_2.$ In fact,
 \[ \langle f_1,B_1\rangle =\langle f,B_1\rangle, \quad {\rm and}\quad \langle f_2,B_2\rangle =\langle f,B_2\rangle,\]
 and hence $ \langle f_1,B_1\rangle B_1+ \langle f_2,B_2\rangle B_2$ is the orthogonal projection of $f$ into the span of $(B_1,B_2),$ and $f_3=f-\langle f_1,B_1\rangle B_1+ \langle f_2,B_2\rangle B_2$ is in the orthogonal complement of the span $(B_1,B_2).$ We therefore have
 \[ \|f\|^2=|\langle f,B_1\rangle|^2+|\langle f,B_2\rangle|^2+\|f_3\|^2\]
 with the consecutive  selections of the parameters $a_1$ and $a_2$ giving rise to  the maximal energy gains in
 $|\langle f,B_1\rangle|^2$ and $|\langle f,B_2\rangle|^2,$ respectively.
We note that $f_3\in I_{\{a_1,a_2\}}.$ That is due to the fact that $f_3\perp {\rm span} \{B_1, B_2\}={\rm span} \{\tilde{k}_{a_1}, \tilde{k}_{a_2}\}.$
 Continuing the above process up to the $n$-th time, we have
 \begin{eqnarray}\label{as defined} f=\sum_{k=1}^n \langle f_k,B_k\rangle B_k + f_{n+1}=\sum_{k=1}^n \langle f,B_k\rangle B_k + f_{n+1},\end{eqnarray}
 where $B_k=B_k^{a_k}, k=1,\cdots n;$ and
 \[ \|f\|^2=\sum_{k=1}^n |\langle f,B_k\rangle|^2 + \|f_{n+1}\|^2.\]
 The feasibility of the maximal selection for each $a_k$ is due to the relation
 \[ \lim_{|b|\to 1-}|\langle f_k, \frac{{e}_b}{\sqrt{1-\sum_{l=1}^{k-1}|\langle e_b,B_l\rangle|^2}} \rangle |=0\]
 proved through the same reasoning as that for the $n=2$ case.  We note that
 $(B_1,\cdots,B_n)$ is a G-S orthonormalization of $(\tilde{k}_{a_1},\cdots,\tilde{k}_{a_n}),$ and
 $f_{n+1}\perp {\rm span}\{B_1,\cdots,B_n\}.$

The above argument is summarized as the Maximal Selection Theorem stated as

\begin{theorem}\label{maximal selection}
For any $f\in \mA$ and positive integer $n$ there exists $a_n\in {\bf D}$ such that
\[ |\langle f,B^{a_n}_{n}\rangle|^2=\sup \{ |\langle f,B_{n}^b\rangle|^2\ :\ b\in {\bf D}\},\]
where for any $b\in \D,$ $(B_1,...,B_{n-1},B^b_n)$ is the G-S orthogonalization of $(\tilde{k}_{a_1},...,\tilde{k}_{a_{n-1}},\tilde{k}_b)$ in accordance with the $n$-tuple $(a_1,\cdots,a_{n-1},b).$
\end{theorem}

We next show

\begin{theorem}\label{convergence} Let $f$ be any function in $\mA.$
Under the maximal selections of the parameters $(a_1,\cdots,a_n,\cdots)$ there holds
\[f=\sum_{k=1}^\infty \langle f,B_k\rangle B_k.\]
\end{theorem}

\noindent{\bf Proof}. We prove the convergence by contradiction. Assume that through a sequence of maximally selected  parameters ${\bf a}=(a_1,\cdots,a_n,\cdots)$ we have
\begin{eqnarray}\label{general} f=\sum_{k=1}^\infty \langle f,B_k\rangle B_k + h, \qquad h\ne 0.\end{eqnarray}
By the Bessel inequality and the Riesz-Fisher Theorem $\sum_{k=1}^\infty \langle f,B_k\rangle B_k \in \mA,$ and hence $h\in \mA.$
We separate the series into two parts
\[ f=\sum_{k=1}^M + \sum_{k=M+1}^\infty \langle f,B_k\rangle B_k + h,\]
where we denote
\[ f_{M+1}=\sum_{k=M+1}^\infty \langle f,B_k\rangle B_k + h =g_{M+1}+h.\]
In the last equality, due to the orthogonality, we may replace $f$ with $f_k$ and have
\[ f_{M+1}=\sum_{k=M+1}^\infty \langle f_k,B_k\rangle B_k + h.\]
Since the set of reproducing kernels of $\mA$ is dense, there exists $a\in \bD$ such that $\delta \triangleq |\langle h,e_a\rangle |>0.$ We can, in particular, choose $a$ to be different from all the selected $a_k$'s in the sequence ${\bf a}.$  The contradiction that we are to explore is in relation to the maximal selections of $a_{M+1}$ when $M$ is large.
Now, on one hand, again by the Bessel inequality, we have
\begin{eqnarray}\label{back} |\langle f_{M+1},B_{M+1}\rangle|\to 0, \qquad {\rm as}\quad M\to 0.\end{eqnarray}

On the other hand we will show, for large $M$,
\begin{eqnarray} \label{backback}|\langle f_{M+1},B_{M+1}^{a}\rangle|> \frac{\delta}{2}.\end{eqnarray}
This is then clearly a contradiction.

The rest part of the proof is devoted to showing (\ref{backback}). Due to the quantity relations
\begin{eqnarray}\label{g+h} |\langle f_{M+1}, B_{M+1}^a\rangle|\ge |\langle h, B_{M+1}^a\rangle|-|\langle g_{M+1}, B_{M+1}^a\rangle|\end{eqnarray} and
\begin{eqnarray}\label{h} |\langle g_{M+1}, B_{M+1}^a\rangle|\leq \|g_{M+1}\|\to 0, \qquad {\rm as}\quad M\to \infty,\end{eqnarray} we see that
if $M$ grows large, the lower bounds of $|\langle f_{M+1}, B_{M+1}^a\rangle|$ depend on those of $|\langle h, B_{M+1}^a\rangle|$. To analyze the lower bounds of  $|\langle h, B_{M+1}^a\rangle|$ we first recall the $h$ is orthogonal to all $B_k.$
 This is seen from the relations
\[ \langle h, B_k\rangle =\lim_{n\to \infty} \langle f-\sum_{l=1}^n \langle f,B_l\rangle B_l, B_k\rangle=0.\]
In the sequel of the proof let the generalized reproducing kernels $\tilde{k}_{a_1}, \cdots, \tilde{k}_{a_n}, \cdots $ be fixed according to the selected sequence ${\bf a}.$
Now for any positive integer $M$ denote by $X_{M+1}^{a}$ the $(M+1)$-dimensional space spanned by $\{{k}_{a}, \tilde{k}_{a_1},\cdots,\tilde{k}_{a_M}\}.$ We have two methods to compute the energy of the projection of $h$ into  $X_{M+1}^{a},$ denoted by $\|h/X_{M+1}^a\|^2.$ One method is based on the orthonormalization $(B_1,\cdots,B_M,B_{M+1}^a).$ In such way, due to the orthogonality of $h$ with $B_1,\cdots,B_M,$ we have
\[ \|h/X_{M+1}^a\|^2=|\langle h,B_{M+1}^a\rangle|^2.\]
The second way is based on the orthonormalization in the changed order $(e_a,\tilde{k}_{a_1},\cdots,\tilde{k}_{a_M}).$ Then we have
\[ \|h/X_{M+1}^a\|^2\ge |\langle h,e_a\rangle|^2=\delta^2.\]
Hence we have, for any $M$, $|\langle h,B_{M+1}^a\rangle|\ge \delta.$ In view of this last estimation and (\ref{h}), (\ref{g+h}), we arrive at the contradiction spelt by (\ref{backback}), (\ref{back}). The proof is thus complete.
\qed

\begin{remark} The above proof is suggested by T. Qian (also see \cite{QianBook1}). \end{remark}

\section {The Rational Orthonormal Systems of $\mA$}

Given any sequence ${\bf a}=(a_1,\cdots,a_n,\cdots)$ in the unit disc ${\bD}$ we correspondingly have the generalized reproducing kernel sequence $(\tilde{k}_{a_1},\cdots,\tilde{k}_{a_n},\cdots),$ and the G-S orthonormalization of the latter, the corresponding orthonormal system $(B_1,\cdots,B_n,\cdots).$ The last system is named as the rational orthonormal system associated with ${\bf a}$, or the Bergman-Takenaka-Malmquist system associated with ${\bf a}.$

We notice that if ${\bf a}$ is not produced under a maximal selection process, then indeed the expansion (\ref{general}) may result in a nontrivial infinite error function $h\ne 0.$
We have shown $h\perp B_k, k=1,\cdots,n. $ Recall that
$$I_{\bf a}=\{f\in \mA : f(\bf a)=\rm0\}.$$
We further show
\begin{theorem}\label{zerobased}
Let ${\bf a}$ be any sequence in $\bD$ and $h$ the infinite error function
\begin{eqnarray}\label{general} h=f-\sum_{k=1}^\infty \langle f,B_k\rangle B_k .\end{eqnarray} Then there holds
$h\in I_{\bf a}.$ Moreover,
 \begin{eqnarray}\label{direct sum}\mA=\overline{{\rm span}}\{B_k\}_{k=1}^\infty \oplus I_{\bf a}.\end{eqnarray} \end{theorem}

\noindent{\bf Proof}. First, the Riesz-Fisher Theorem gives
$$ h=f-\sum_{k=1}^\infty \langle f, B_k\rangle B_k\in \mA.$$ Next, we show that $h$ has ${\bf a}$ as its zero set.
The fact $h\perp B_k$ for all $k$ implies $h\perp \tilde{k}_{a_k}$ for all $k.$ This implies that $h$ itself and all its derivatives up to the order $l(a_k)-1$ have $a_k$ as zero. So, $h$ has at least $l(a_k)$-order zero $a_k.$ This is for all positive integers $k.$ Note that the proof does not exclude the case $h=0.$  As a consequence of the orthogonality the LHS of (\ref{direct sum}) is a subset of the RHD. The inverse inclusion is obvious. \qed\\

As a consequence we have
\begin{coro}\label{basis}
 $\{B_k\}_{k=1}^\infty$ is a basis if and only if $I_{\bf a}=\{0\}.$
\end{coro}

The following result is contained in \cite{Duren}.

\begin{lemma}\label{5.3}
If ${\bf a}$ is the zero set of some non-trivial function $f\in \mA,$ then
$\sum_{k=1}^\infty (1-|a_k|)^2<\infty,$ and the H-product $H_{\bf a}$ is well defined, and
$$ \|f/H_{\bf a}\|\leq C \|f\|.$$
As a consequence, $I_{\bf a}\subset H_{\bf a}\mA.$
\end{lemma}

 It is an interesting question that for what ${\bf a}$ there holds  $I_{\bf a}\ne \{0\}.$ This has been an open question till now. Some results on this may be found in \cite{B}. Studies show that conditions only on the magnitudes of the $a_k$'s are not sufficient to guarantee ${\bf a}$ to be the zero set of some non-trivial $f\in \mA.$ Some results along this line can be read off from the literature. Before we state our results we introduce a probability on an infinite product space made from the arguments of the $a_j$'s \cite{prob, QWang}.\\

\noindent{\bf Definition} Let $\Omega=\prod_{j=1}^\infty [0,2\pi ).$ Let $\mu_j$ be the normalized Lebesgue measure on the $j$-th factor of $\Omega.$ Then let $\mu$ be the probability measure on $\Omega$ such that $\mu=\prod_{k=1}^\infty \mu_j.$ \\

 Let $\{r_j\}_{j=1}^\infty$ be an increasingly ordered sequence (allowing repetition) in $[0,1)$ satisfying
  \begin{eqnarray*}\label{hyper}\sum_{j=1}^\infty (1-r_j)^2<0.\end{eqnarray*}

  Consider the map of $\Omega$ into the holomorphic function set $H(\bD)$ defined by $\omega \to H_{\bf a}$ such that
$a_j=r_je^{i\omega_j}$ and $\omega =\{\omega_j\}_{j=1}^\infty.$ Then under the condition
\begin{eqnarray}\label{epsilon} \lim \sup_{\epsilon \to 0+} \frac{\sum_{j=1}^\infty (1-r_j)^{1+\epsilon}}{\log\frac{1}{\epsilon}}<1/4,\end{eqnarray}
there holds that $H_{\bf a}\in \mA$ for $\mu$-a.e. $\omega$ \cite{B,prob}.

Based on this result, as well as Corollary \ref{basis} and Lemma \ref{5.3}, we conclude
\begin{theorem}
Under the condition (\ref{epsilon}) $\{B_j\}_{j=1}^\infty$ is almost surely not a basis.
\end{theorem}

\section {Convergence Rate of Pre-Orthogonal Adaptive Fourier Decomposition}

 The POAFD approach is not classified into any existing category of the greedy algorithms. The differences between the AFD methods (including POAFD) and the greedy algorithms include (i) in the AFD methods the parameters can be selected repeatedly; and (ii) thus at each of the iterative steps the maximal projection can be attained. It, therefore, offers approximations better than greedy algorithms \cite{Q2}. Below we prove the corresponding convergence rate.

 Before we prove the convergence rate estimation we recall the following lemma whose proof can be found in the greedy algorithm literature.

\begin{lemma}\label{lemma} Assume that a sequence of positive numbers satisfies the conditions
\[ d_1\leq A, \quad d_{k+1}\leq d_k(1-\frac{d_k}{r_k^2A}), \qquad k=1,2,\cdots\]
Then there holds
\[ d_k\leq \frac{A}{1+\sum_{l=1}^k \frac{1}{r_l^2}}, \qquad k=1,2,\cdots\]
\end{lemma}

For $M>0$ we will be working with the subclass ${\mA}_M$ of $\mA$ defined as
\[ {\mA}_M=\{ f\in \mA\ :\ \exists \{b_1,\cdots,b_n,\cdots\}, f=\sum_{l=1}^\infty c_l{e}_{b_l}, \sum_{l=1}^\infty |c_l|\leq M\}.\]

\begin{theorem} Let $f=\sum_{l=1}^\infty c_l{e}_{b_l} \in {\mA}_M, \sum_{l=1}^\infty |c_l|\leq M,$ and $f_n$ be the orthogonal standard remainder as defined in \rm{(\ref{as defined})} corresponding to the maximal selections of the $a_k$'s, then there exists estimation
 \[ \|f_k\|\leq M\left(1+\sum_{l=1}^k\left(\frac{1}{r_l}\right)^2\right)^{-\frac{1}{2}},\]
 where $r_k=\sup \{ r_k(b_l)\ : \ l=1,2,\cdots \},$ and $r_k(b_l)=\sqrt{1-\sum_{t=1}^{k-1}|\langle {e}_{b_l},B_t\rangle|^2}=\|e_{b_l}-\sum_{t=1}^{k-1}\langle {e}_{b_l},B_t\rangle B_t\|.$
\end{theorem}

\noindent{\bf Proof}. We start from the inequality chain
\begin{eqnarray*}
|\langle f_k,B_k\rangle|&\ge& \sup \{ |\langle f_k,B_k^a\rangle|\ :\ a\in {\bf D}\}\\
&\ge & \sup \{ |\langle f_k,B_k^{b_l}\rangle|\ :\ l=1,2,\cdots\}\\
&= & \sup \{ \frac{|\langle f_k,e_{b_l}\rangle|}{r_k(b_l)}\ :\ l=1,2,\cdots\}\\
&\ge & \frac{1}{r_k} \sup \{{|\langle f_k,e_{b_l}\rangle|}\ :\ l=1,2,\cdots\}\\
&\ge & \frac{1}{r_k M} |\langle f_k,\sum_{l=1}^\infty  c_l{e}_{b_l}\rangle|\\
&= & \frac{1}{r_k M} |\langle f_k,f\rangle|\\
&= & \frac{1}{r_k M} |\langle f_k,f_k\rangle| \qquad ({\rm since}\ f_k\perp (f-f_k))\\
&= & \frac{\|f_k\|^2}{r_k M}.
\end{eqnarray*}
Substituting this inequality into the relation
\[ \|f_{k+1}\|^2=\|f_k\|^2-|\langle f_k,B_k\rangle|^2,\]
we have
\[ \|f_{k+1}\|^2\leq \|f_k\|^2\left(1-\frac{1}{(r_kM)^2}\|f_k\|^2\right).\]

By invoking Lemma \ref{lemma}, we have the desired estimation
\[ \|f_k\|\leq \frac{M^2}{\sqrt{1+\sum_{l=1}^k\frac{1}{r_l^2}}}.\]
The proof is complete. \qed\\

\noindent{\bf Remark} Since $0<r_k\leq 1,$ we have, at least
\[ \|f_k\|\leq \frac{M^2}{\sqrt{k}}.\]
Coincidentally this is the same bound as for the Shannon expansion. However, Shannon expansion treats bandlimited entire functions with great smoothness that is what is usually needed for good convergence rates.
On the other hand, the Bergman space contains functions that blow up at the boundary.

\section {The Bergman Space on the Upper Half Complex Plane}
The theory on the upper half complex plane $\bf{C}_+$ is a close analogy with the one for the unit disc.  Denote by $\mathbb{A}^2({\bf{C}_+})$ the square integrable Bergman space of the upper half complex plane $\bf{C}_+$, that is,
$$
 \mathbb{A}^2({\bf{C}_+})=\{f: {\bf{C}_+}\to {\bf C}\ | \ f\ {\rm is \ holomorphic \ in \ \bf{C}_+, \ and\ }\ \| f\|_{\mA({\bf{C}_+})}^2=\int_{\bf{C}_+} |f(z)|^2dA<\infty\},
 $$
The reproducing kernel $k^+_a$ of $\mA({\bf{C}_+})$ at the point $a\in \bf{C}_+$ is given by $ k^+_a(z)=\frac{-1}{(z-\overline{a})^2}$ \cite{Chinese,E}. Moreover, $ \| k^+_a\|^2= k^+_a(a)=\frac{1}{(2\rm{Im}\it a)^{\rm{2}}}$ and the normalized kernel $e^+_{a}=\frac{k^+_a(a)}{\| k^+_a\|}=\frac{-2\rm{Im}\it a}{(z-\overline{a})^2}.$ In the later part of this section we will write, for the simplicity but with an abuse of notation, $k^+_a, e^+_a,$ etc., as $k_a, e_a,$ etc.

\begin{lemma}\label{lemmaC}
$\Span\{\frac{1}{(z-\overline{a})^2}\ |\ a\in \bf{C}_+\}$ is dense in $\mA({\bf{C}_+}).$
\end{lemma}
\noindent{\bf Proof}. Let $\mathcal{A}=\overline{\Span}\{\frac{1}{(z-\overline{a})^2}\ |\ a\in \bf{C}_+\},$ we claim that $\mathcal{A}=\mathbb{A}^2({\bf{C}_+}).$ If the last relation did not hold, then $$\mathbb{A}^2({\bf{C}_+})=\mathcal{A}\oplus\mathcal{A^{\perp}} \ \ \rm{and}\ \  \mathcal{A^{\perp}}\neq \{0\}.$$
Thus, there exists $f\in \mathcal{A^{\perp}}$ and $f\neq 0.$ In such case due to the reproducing kernel property we have $<f,\frac{1}{(z-\overline{a})^2}>=f(a)=0$ for every $a\in \bf{C}_+.$ It follows that $f= 0.$
It is a contradiction. Therefore, $\overline{\Span}\{\frac{1}{(z-\overline{a})^2}\ |\ a\in \bf{C}_+\}=\mathbb{A}^2({\bf{C}_+}).$ \qed\\ \par

\begin{lemma}\label{premaximal selectionC}
For any $f\in \mA({\bf{C}_+})$ there exists $a\in {\bf {C}_+}$ such that
\[ |\langle f,e_{a}\rangle|^2=\max \{ |\langle f,e_{b}\rangle|^2\ :\ b\in {\bf {C}_+}\}.\]
\end{lemma}

\noindent{\bf Proof}
The Cauchy-Schwarz inequality gives
\[|\langle f,e_b\rangle|\leq \|f\|.\]
Therefore $ |\langle f,e_b\rangle| $ has a finite upper bound. In particular,
$\sup \{|\langle f, e_b\rangle|\ |\ b\in {\bf {C}_+}\}\leq \|f\|.$ There exists a sequence $\{b_k\}_{k=1}^\infty$ such that
\begin{eqnarray}\label{due to} \lim_{k\to \infty}|\langle f, e_{b_k}\rangle|=\sup \{|\langle f, e_b\rangle|\ |\ b\in {\bf {C}_+}\}.\end{eqnarray}
We will prove there exists a subsequence $\{b_{k_l}\}_{l=1}^\infty$ converging to $\tilde{b}\in {\bf {C}_+}.$
 Then we can conclude that $|\langle f, e_{\tilde{b}}\rangle|=\lim_{l\to \infty}|\langle f, e_{{b_{k_l}}}\rangle|=\sup \{|\langle f, e_b\rangle|\ |\ b\in {\bf {C}_+}\},$ as desired.

The strategy is to show:

For any $\epsilon>0$ there exists an open neighborhood $B(\delta, R)$ of the boundary of ${\bf C}^+,$ where $B(\delta, R)=\{b\in {\bf C}^+\ |\ {\rm Im}b <\delta\}\cup \{b\in {\bf C}^+\ |\ |b|>R\},$ such that $|\langle f, e_b\rangle|<\epsilon$ whenever $b\in B(\delta, R).$ As a consequence of this, as well as of the relation (\ref{due to}), when $k$ grows large those $b_k$ must stay in a compact subset of ${\bf C}^+,$ and thus one can choose a subsequence of $\{b_{k_l}\},$ according to the Bolzano-Weierstrass Theorem, converging to a point $\tilde{b}$ in ${\bf C}^+.$

For $\epsilon >0,$ since $\Span\{\frac{1}{(z-\overline{a})^2}\ |\ a\in \bf{C}_+\}$ is dense in $\mA({\bf{C}_+})$, there exists
\[g=\sum_{k=1}^{m}\frac{c_k}{(z-\overline{a}_k)^2}\]
such that
\[ \|f-g\|\leq \epsilon /2.\]
Denote min$\{\rm{Im}\it a_k\}_{k=1}^{m}=\delta_1>0.$ As before, the triangle inequality and the Cauchy-Schwarz inequality give
\[ |\langle f,e_b\rangle|\leq |\langle g,e_b\rangle|+\epsilon/2.\] \par
On one hand, there exists $\delta >0$ such that
${\rm Im} b <\delta $ implies $|\langle g,e_b\rangle|<\epsilon/2.$ In fact, we have
$$
|\langle g,e_b\rangle|=2\rm Im\it b\ |\sum_{k=1}^{m}\frac{c_k}{(b-\overline{a}_k)^{\rm 2}}|\leq \rm2\rm Im\it b\ \sum_{k=1}^{m}\frac{|c_k|}{|b-\overline{a}_k|^{\rm 2}}\leq \rm2\rm \delta \ \sum_{k=1}^{m}\frac{|c_k|}{(\frac{\delta_1}{\rm2})^{\rm2}}<\epsilon/2,
$$
if $\delta $ is small enough.

On the other hand, let $\max \{|a_k|\}_{k=1}^{m}=R_1>0,$ We will show there exists $R$ such that $|b|>R$ implies $|\langle g,e_b\rangle|<\epsilon/2.$ \\
Let $|b|>4R_1.$ There holds $|b-\overline{a}_k|\ge |b|/2.$ Hence
$$
|\langle g,e_b\rangle|\leq \rm2\rm Im\it b\ \sum_{k=\rm1}^{m}\frac{|c_k|}{|b-\overline{a}_k|^{\rm {\rm 2}}}\leq \rm 2|b|\ \sum_{\it k=\rm1}^{\it m}\frac{4|\it c_k|}{|b|^2}=\frac{1}{|b|}\sum_{\it k=\rm 1}^{\it m}8|\it c_k|.
$$
So, if $|b|>R,$ and $R$ is large enough and larger than $4R_1,$ then the last quantity is dominated by $\epsilon/2.$
\qed \par
 Note that although we did not define $\partial {\bf C}^+$ but we define the open neighborhoods of $\partial {\bf C}^+,$ viz. $B(\delta,R).$ By $b$ being further close to $\partial {\bf C}^+$ we mean that $b$ is in $B(\delta',R')$ with $\delta'<\delta,\ R'>R.$ The above proved is regarded as uniform boundary vanishing property of $\langle f, e_b\rangle,$  and denoted  $$\lim_{b\to\partial\bf{C}_+}|\langle f,e_b\rangle|=0.$$ This concludes that $|\langle f,e_b\rangle|$ can attain the maximum value at an interior point (the maximal selection principle for $B_1=e_{a_1}$ in the upper-half plane). Next we will develop the POAFD algorithm together with the adaptive TM system of the Bergman space in the upper-half complex plane.

For $f_1\in \mA({\bf{C}_+})$ choose $a_1$ such that
\[  |\langle f_1,e_{a_1}\rangle|^2=\max \{ |\langle f_1,e_{b}\rangle|^2\ :\ b\in {\bf {C}_+}\}.\]
Fixing $a_1,$ we have, by denoting $e_{a_1}=B_1,$
 \begin{eqnarray}\label{f3C} f(z)=\langle f_1, B_1\rangle B_1(z) + \langle f_2, B_2^b\rangle B_2^b(z) + f_3,\end{eqnarray}
 where $(B_1,B_2^b)$ is a G-S orthonormalization of $(B_1, e_{b}),$ where $b\ne a_1.$ A close look at $\langle f_2, B_2^b\rangle$ gives

 \begin{eqnarray*}
 \langle f_2,B_2^b\rangle = \langle f_2, \frac{e_b}{\sqrt{1-|\langle e_b,B_1\rangle|^2}} \rangle=\langle f_2, \frac{e_b}{\sqrt{1-|\langle e_b,B_1\rangle|^2}}- e_b \rangle + \langle f_2, e_b \rangle= I_1(b)+I_2(b)
 \end{eqnarray*}
Being similar with the unit disc case we have
 \[ |I_1(b)|\leq \int_{\bf{C}_+} |f_2(z)||e_b(z)|\left(\frac{1}{\sqrt{1-|\langle e_b,B_1\rangle|^2}}-1\right)dA(z).\]

 By invoking the reproducing kernel property, if $\rm{Im}\it b$ is close to 0, then
 \[|\langle e_b,B_1\rangle|=\frac{(2\rm{Im}\it {a_{\rm 1}})(\rm{2Im}\it b)}{|a_1-\overline{b}|^2}\leq \epsilon; \]
  and if $|b|$ is close to $\infty$, then
 \[|\langle e_b,B_1\rangle|\leq \frac{\rm Im \it a_{\rm 1}}{|b|}\leq \epsilon.\]
 Therefore,
 \[ \lim_{b\to \partial {\bf C}_+} I_1(b)=0.\]
The maximal selection principle implies $\lim_{b\to \partial {\bf C}_+} I_2(b)=0.$ To summarize, we have
 \begin{eqnarray}\label{same reasoningC} \lim_{b\to \bf \partial C_+} |\langle f_2,B_2^b\rangle|=0.\end{eqnarray}

 Hence there exists $a_2\in \bf{C}_+$ such that
 \[ |\langle f_2,B_2^{a_2}\rangle|=\max \{ |\langle f_2,B_2^b\rangle|\ :\ b\in \bf{C}_+\}.\]

There also may happen on the upper-half complex plane, through a limiting sequence $b_k$ gaining the maximal value, $\lim b_k=a_1.$ In such case, we define $B_2=B_2^{a_2}=\lim B_2^{b_k},$ that involves a directional derivative of the reproducing kernel $k_{a_1}$ as explained in \S 3 and \S 4.  After $B_1, B_2$ being well defined, we have, likewise,
\[ \langle f_1,B_1\rangle =\langle f,B_1\rangle, \quad {\rm and}\quad \langle f_2,B_2\rangle =\langle f,B_2\rangle,\]
 and
 $f_3\perp {\rm span} \{B_1, B_2\}={\rm span} \{\tilde{k}_{a_1}, \tilde{k}_{a_2}\}.$
 We therefore have
 \[ \|f\|^2=|\langle f,B_1\rangle|^2+|\langle f,B_2\rangle|^2+\|f_3\|^2.\]
 Repeating this process up to the $n$-th time, we produce a decomposition result essentially in agreement with the unit disk case:
 \begin{eqnarray}\label{Was defined} f=\sum_{k=1}^n \langle f_k,B_k\rangle B_k + f_{n+1}=\sum_{k=1}^n \langle f,B_k\rangle B_k + f_{n+1},\end{eqnarray}
 where $B_k=B_k^{a_k}, k=1,\cdots n;$ and
 \[ \|f\|^2=\sum_{k=1}^n |\langle f,B_k\rangle|^2 + \|f_{n+1}\|^2.\]
 The feasibility of the maximal selection for each $a_k$ on $\bf C_+$ is due to the relation
 \[ \lim_{b\to \bf \partial C_+}|\langle f_k, \frac{{e}_b}{\sqrt{1-|\sum_{l=1}^{k-1}\langle e_b,B_l\rangle|^2}} \rangle |=0\]
 proved through similar estimation as in the unit disc.  We note that
 $(B_1,\cdots,B_n)$ is the G-S orthonormalization of $(\tilde{k}_{a_1},\cdots,\tilde{k}_{a_n}),$ and
 $f_{n+1}\perp {\rm span}\{B_1,\cdots,B_n\}.$\par
According to the above argument, we have a maximal selection theorem on $\bf{C}_+$ analogous with the one on $\D:$

\begin{theorem}\label{maximal selectionC}
For any $f\in \mA({\bf{C}_+})$ and positive integer $k$ there exists $a_k\in {\bf {C}_+}$ such that
\[ |\langle f,B^{a_k}_{k}\rangle|^2=\sup \{ |\langle f,B_{k}^b\rangle|^2\ :\ b\in {\bf {C}_+}\},\]
where for any $b\in \bf{C}_+,$ $(B_1,...,B_{k-1},B^b_k)$ is the G-S orthogonalization of $(B_1,...,B_{k-1},\tilde{k}_b)$ in accordance with the $k$-tuple $(a_1,\cdots,a_{k-1},b).$
\end{theorem}

With a similar proof as for Theorem \ref{convergence} we have

\begin{theorem}\label{convergenceC} Let $f$ be any function in $\mA({\bf{C}_+}).$
Under the maximal selections of the parameters $(a_1,\cdots,a_n,\cdots)$ there holds
\[f=\sum_{k=1}^\infty \langle f,B_k\rangle B_k.\]
\end{theorem}

The above convergence theorem has a similar proof as in the unit disk case.

\section{Weighted Bergman Spaces $\mA_\alpha$ with $-1<\alpha<\infty$}
In this section we study general weighted Bergman spaces $\mA_\alpha\ ,\ -1<\alpha<\infty.$ We adopt the notation
$$ \mathbb{A}_\alpha^2({\bf D})=\{f: {\D}\to {\bf C}\ | \ f\ {\rm is \ holomorphic \ in \ \D, \ \rm and\ }\ \| f\|_{\mA_\alpha({\bf D})}^2=\int_{\bf D} |f(z)|^2dA_\alpha<\infty\},$$
where $dA_\alpha(z)=(1+\alpha)(1-|z|^2)^{\alpha}dA(z).$ With the norm $\| \cdot\|_{\mA_\alpha({\bf D})}^2$ the space $\mA_\alpha$ is a reproducing kernel Hilbert space \cite{Zhu}. The reproducing kernel and its norm are given, respectively, by
$$ k_a^{\alpha}(z)=\frac{1}{(1-\overline{a}z)^{2+\alpha}}\ \  \rm and \ \ \it \| k_a^{\alpha}\|^{\rm 2}= k_a^{\alpha}(a)=\frac{\rm 1}{(\rm 1-|\it a|^{\rm 2})^{\rm 2+\alpha}}.$$
In this section we again denote by $(B_1,...,B_k)$ the G-S orthonormalization of $(\tilde{k}^\alpha_{a_1},...,\tilde{k}^\alpha_{a_k}),$ where $\tilde{k}^{\alpha}_{a_l}$ is defined as in (\ref{derivative}) depending on the multiple of $a_l$ in
$(a_1,...,a_l), 1\leq l\leq k.$

 If  $f(z)=\sum_{k=0}^{\infty}a_k z^k\in \mA_\alpha, $ then simple computation gives
\begin{equation}\label{gamma}
\| f\|_{\mA_\alpha}^2=\sum_{k=0}^{\infty}\frac{k!\Gamma(\alpha+2)}{\Gamma(k+\alpha+2)}|a_k|^2=\sum_{k=0}^{\infty}\frac{k!}{(\alpha+1+k)(\alpha+k)\cdots(\alpha+2)}|a_k|^2.
\end{equation}
The above energy formula is in terms of the multiplier $\frac{k!\Gamma(\alpha+2)}{\Gamma(k+\alpha+2)}.$ When $\alpha$ is an integer, it is\par
\noindent$\sum_{k=0}^{\infty}\frac{\Gamma(\alpha+2)}{(k+\alpha+1)\cdots(k+1)}|a_k|^2$
and, when $\alpha=0,$ reduces to $\sum_{k=0}^{\infty}\frac{1}{k+1}|a_k|^2,$ varifying the case of the classical Bergman space.
It can be easily shown that the maximal selection principle holds for $\mA_\alpha$:

\begin{theorem}\label{weighted bergman msp}
For any $f\in \mA_{\alpha}$ and positive integer $k$ there exists $a_k\in {\bf D}$ such that
\[ |\langle f,B^{a_k}_{k}\rangle|^2=\sup \{ |\langle f,B_{k}^b\rangle|^2\ :\ b\in {\bf D}\},\]
where for any $b\in \D,$ $(B_1,...,B_{k-1},B^b_k)$ is the G-S orthonormalization of $(B_1,...,B_{k-1},\tilde{k}^{\alpha}_b)$ being in accordance with the $k$-tuple $(a_1,\cdots,a_{k-1},b).$
\end{theorem}

The proof is similar to that of Theorem \ref{maximal selection}. As a consequence, POAFD can be applied to $\mA_\alpha$ as for the case $\alpha=0$ developed in the former sections.\par
The proof of Theorem \ref{convergence} can be adapted to the weighted Bergman space cases. There holds
\begin{theorem}\label{weighted bergman convergence}
Let $f$ be any function in $\mA_{\alpha}.$
Under the maximal selections of the parameters $(a_1,\cdots,a_n,\cdots)$ there holds
\[f=\sum_{k=1}^\infty \langle f,B_k\rangle B_k.\]
\end{theorem}

Next we will give a slight discuss on the spaces $\mA_{\alpha},\ \ \alpha >-1.$
\begin{lemma}\label{weighted bergman series equivalance}
$\sum_{k=0}^{\infty}\frac{k!\Gamma(\alpha+2)}{\Gamma(\alpha+2+k)}|a_k|^2$ is equivalent with $\sum_{k=0}^{\infty}\frac{1}{(k+1)^{\alpha+1}}|a_k|^2,$ that is, the two series converge or diverge simultaneously.
\end{lemma}

\noindent{\bf Proof}. Let $\sum_{k=0}^{\infty}\frac{k!\Gamma(\alpha+2)}{\Gamma(\alpha+2+k)}=\sum u_k$ and $\sum_{k=0}^{\infty}\frac{1}{(k+1)^{\alpha+1}}=\sum v_k. $ Then $\sum u_k$ and $\sum v_k$ are series of positive terms. We will show that \begin{eqnarray}\label{lemma basic}
\lim_{k\to \infty}\frac{u_k}{v_k}=C \in (0,\infty).
\end{eqnarray}
 Based on (\ref{lemma basic}) we can claim the desired result through the comparison principle for series of positive terms. We are to show (\ref{lemma basic}).
In fact, by the Stirling formula
\[ \Gamma (z)=\sqrt{\frac{2\pi}{z}}\left(\frac{z}{e}\right)^z[1+O(\frac{1}{z})],\] we have
\begin{eqnarray*}
\lim_{k\to \infty}\frac{u_k}{v_k} &=& \Gamma(\alpha+2)\ \lim_{k\to \infty}\frac{\Gamma(k+1)(k+1)^{\alpha+1}}{\Gamma(\alpha+2+k)}\\
&=& \Gamma(\alpha+2)\lim_{k\to \infty}\ \frac{\sqrt{\frac{2\pi}{k+1}}(\frac{k+1}{e})^{k+1}[1+O(\frac{1}{k+1})](k+1)^{\alpha+1}}{\sqrt{\frac{2\pi}{\alpha+2+k}}(\frac{\alpha+2+k}{e})^{\alpha+2+k}[1+O(\frac{1}{\alpha+2+k})]}\\
&=& \Gamma(\alpha+2)\lim_{k\to \infty}\sqrt{1+\frac{\alpha+1}{k+1}}e^{\alpha+1}(1-\frac{\alpha+1}{\alpha+2+k})^{\frac{\alpha+2+k}{\alpha+1}(\alpha+1)}\frac{1+O(\frac{1}{k+1})}{1+O(\frac{1}{\alpha+2+k})}\\
&=& \Gamma(\alpha+2).  \qed\\
\end{eqnarray*}

\begin{theorem}\label{weighted bergman new}
Let $-1<\alpha<\infty,$ then \par
(i) $\mA_\alpha$ strictly increases as $\alpha$ increases, i.e. $-1<\alpha_1<\alpha_2<\infty$ implies $\mA_{\alpha_1} \subsetneqq \mA_{\alpha_2}.$ \par
(ii) $H^2(\D)\subsetneqq \bigcap_{\alpha>-1} \mA_{\alpha}.$ \par
(iii) $\bigcup_{-1<\gamma <\alpha}\mA_{\gamma}\subset \mA_{\alpha}\subsetneqq \bigcap_{\beta \ge \alpha}\mA_\beta.$\par
(iv) For $-1<\alpha<\beta<\infty,\ P(\D),$ $H^2(\D)$ and $\mA_\alpha$ are all dense in $\mA_\beta$ (in the norm of $\mA_\beta$).
\end{theorem}

\noindent{\bf Proof}. (i) First, supposed $f=\sum_{k=0}^{\infty}a_kz^k \in \mA_{\alpha_1},$ according to lemma \ref{weighted bergman series equivalance}, we have
\[\sum_{k=0}^{\infty}\frac{|a_k|^2}{(k+1)^{\alpha_2+1}}<\sum_{k=0}^{\infty}\frac{|a_k|^2}{(k+1)^{\alpha_1+1}}<\infty.\]
Then $f\in \mA_{\alpha_2}.$ Secondly, let $f=\sum_{k=0}^{\infty}a_kz^k \in \mA_{\alpha_2},$ where $|a_k|^2=\frac{1}{(k+1)^{1+\delta}},\ -1-\alpha_2<\delta<-1-\alpha_1.$ Taking $\delta =-1-\frac{\alpha_1+\alpha_2}{2},$ an easy computation gives
\[ \sum_{k=0}^{\infty} \frac{1}{(k+1)^{\alpha_2+1}}|a_k|^2=\sum_{k=0}^{\infty} \frac{1}{(k+1)^{1+\frac{\alpha_2-\alpha_1}{2}}}<\infty,\]
while
\[ \sum_{k=0}^{\infty} \frac{1}{(k+1)^{\alpha_1+1}}|a_k|^2=\sum_{k=0}^{\infty} \frac{1}{(k+1)^{1+\frac{\alpha_1-\alpha_2}{2}}}=\infty.\]
This shows that $f(z)=\sum_{k=0}^\infty a_kz^k\in \mA_{\alpha_2}\setminus \mA_{\alpha_1}.$\\

(ii) Supposed $f=\sum_{k=0}^{\infty}a_kz^k \in H^2(\D),$ then
\begin{eqnarray*}
\sum_{k=0}^{\infty}\frac{k!\Gamma(\alpha+2)}{\Gamma(\alpha+2+k)}|a_k|^2 &=& \sum_{k=0}^{\infty}\frac{k!}{(\alpha+1+k)(\alpha+1+k-1)\cdots(\alpha+2)}|a_k|^2\\
&=& \sum_{k=0}^{\infty}\frac{k}{\alpha+1+k}\frac{k-1}{\alpha+1+k-1}\cdots \frac{1}{\alpha+2}|a_k|^2\\
&<&\sum_{k=0}^{\infty}|a_k|^2\\
&<&\infty,
\end{eqnarray*}
showing that $f \in \mA_\alpha.$ We further show that there exists $g \in \mA_\alpha$ for all $\alpha>-1,$ but $g \notin H^2(\D).$ Let $g=\sum_{k=0}^\infty b_kz^k,$ where $|b_k|^2=\frac{1}{k+1},$ then $\sum_{k=0}^{\infty}\frac{1}{(k+1)^{\alpha+1}}|b_k|^2<\infty.$ But $\sum_{k=0}^{\infty}|b_k|^2=\sum_{k=0}^{\infty}\frac{1}{k+1}=\infty,$ so $g \notin H^2(\D).$ \\
(iii) Due to the strict increasing property along with the parameter, the inclusion relations of the spaces are obvious. To show the non-identical relation,
 we take $|a_k|^2=(k+1)^\alpha.$ Then
\[ \sum_{k=0}^\infty \frac{1}{(k+1)^{\alpha+1}}|a_k|^2=\infty;\]
while
\[ \sum_{k=0}^\infty \frac{1}{(k+1)^{\alpha+\delta+1}}|a_k|^2<\infty, \quad \beta-\alpha =\delta>0.\]
Those together mean that $f(z)=\sum_{k=0}^\infty a_kz^k\in \bigcap_{\beta \ge \alpha}\mA_\beta\setminus \mA_\alpha.$ \\
(iv) Since the set of polynomials is dense in $\mA_\beta,$ all function sets between the set of polynomials and the $\mA_\beta$ are dense in the latter.   \qed\\

Next we study a sequence of unbounded holomorphic functions that belong to different weighted Bergman spaces. They are by themselves interesting, as well as useful.

\begin{theorem}\label{beta function}
$f_{\beta}(z)=\frac{1}{(1-z)^{2+\beta}},$ when\\
(i) $\beta>-\frac{3}{2},$ $f_{\beta}(z) \in \mA_{\alpha}\setminus \mA_{2+2\beta},$ $\alpha>2+2\beta.$\\
(ii) $\beta=-\frac{3}{2},$ $f_{\beta}(z) \in \mA_{\alpha}\setminus H^2,$ for any $\alpha>-1.$\\
(iii) $\beta<-\frac{3}{2},$ $f_{\beta}(z) \in H^2.$
\end{theorem}

\noindent{\bf Proof}. We first consider the case $\beta>-\frac{3}{2}.$
\begin{eqnarray*}\label{math analysis}
\int_{\D}\left|\frac{1}{(1-z)^{2+\beta}}\right|^2dA_{\alpha} &=& \int_{\D}\left(\frac{1}{|1-z|^2}\right)^{2+\beta}dA_{\alpha}\notag \\
&=& \int_0^1 r(1+\alpha)(1-r^2)^{\alpha}dr \int_{-\pi}^{\pi}\frac{1}{(1-2r\cos t+r^2)^{2+\beta}}dt\notag \\
&=& 2\int_0^1 r(1+\alpha)(1-r^2)^{\alpha}G(r)dr\\
&=& 2\left(\int_0^{\frac{1}{2}}+\int_{\frac{1}{2}}^1\right) r(1+\alpha)(1-r^2)^{\alpha}G(r)dr\\
&=& I_1+I_2,
\end{eqnarray*}
where
\[ G(r)=\int_0^\pi \frac{1}{(1-2rcost+r^2)^{2+\beta}}dt,\]
and $I_1$ and $I_2$ correspond to the integrals for $r$ in, respectively, $[0,1/2]$ and $[1/2,1].$
First consider $I_1.$ Since $2+\beta>0,$ $1-2r\cos t+r^2\ge (1-r)^2\ge 1/4$ and $\alpha >-1,$
the integral $I_1$ is finite. We are reduced to estimate $I_2.$ For this purpose we separate the inner integral
$G(r)$ into two parts:
\begin{eqnarray*}
G(r)&=&\int_0^\delta \frac{1}{(1-2r\cos t+r^2)^{2+\beta}}dt+\int_\delta^\pi \frac{1}{(1-2r\cos t+r^2)^{2+\beta}}dt\\
&=& J_1(r)+J_2(r),\end{eqnarray*}
where $\delta$ is chosen according to the relation $\sin\frac{\delta}{2}=1-r, 0<\delta<\pi.$ In the sequel $C$ represents a constant that can be different from time to time.
Due to the relation $2+\beta>0$ and $1-2r\cos t+r^2=(1-r)^2+4r\sin^2\frac{t}{2}\leq 5(1-r)^2,$ $0\leq t\leq \delta,$ we have
$$ J_1(r)\ge  C\frac{\delta}{ (1-r)^{4+2\beta}}\ge \frac{C}{(1-r)^{3+2\beta}}.$$
On the other hand,
since $1-2r\cos t+r^2=(1-r)^2+4r\sin^2\left(\frac{t}{2}\right)\ge (1-r)^2, 4+2\beta >-1,$  we have
\begin{eqnarray*}
J_1(r)&\leq& \int_0^\delta \frac{C dt}{(1-r)^{4+2\beta}}\\
&\leq& \frac{C}{(1-r)^{3+2\beta}}.\end{eqnarray*}
For $J_2(r),$ we use the estimate  $1-2r\cos t+r^2=(1-r)^2+4r\sin^2\left(\frac{t}{2}\right)\ge 2\sin^2\left(\frac{t}{2}\right)$ and thus have
\begin{eqnarray*}
J_2(r)&\leq& \int_\delta^\pi \frac{C}{\sin^{4+2\beta}\left(\frac{t}{2}\right)}dt\\
&\leq& C \left[\frac{1}{t^{3+2\beta}}\right]_\delta^\pi\\
&\leq& \frac{C}{(1-r)^{3+2\beta}}.
\end{eqnarray*}
The above estimates for $J_1(r)$ and $J_2(r)$ show that $G(r)$ has singularity only at $r$ being close to $1,$ and
the quantity of the singularity is exactly $(1-r)^{3+2\beta}.$ Therefore, if and only if $\alpha>2+2\beta$ we have
\[ \int_0^\pi r(1-r^2)^\alpha G(r)dr<\infty.\]
This shows that $f_\beta(z)\in \mA_\alpha\setminus \mA_{2+2\beta}$ for $\alpha>2+2\beta.$ The discussion for the case (i) is complete.\\

Now we consider the case (ii) $\beta =-\frac{3}{2}.$ Let $\alpha$ be any real number larger than $-1.$ We are to show $f_\beta\in \mA_\alpha.$ For such $\alpha$ there holds $\alpha =-1+\delta,\ \delta>0.$ Then $\beta_1=-\frac{3}{2}+\frac{\delta}{4}$ belongs to the first case (i). Therefore, for any $\alpha_1> 2+2\beta_1=
2+2[-\frac{3}{2}+\frac{\delta}{4}]=-1+\delta/2,$ the function $f_{\beta_1}\in \mA_{\alpha_1}.$ Find small enough $\alpha_1$ such that $\alpha_1<\alpha=-1+\delta.$ Then we have  $f_{\beta_1}\in \mA_{\alpha_1}\subset \mA_\alpha.$ Since $|f_\beta|\leq |f_{\beta_1}|,$ we also have $f_\beta\in \mA_\alpha.$ Next we need to show that for $\beta=-\frac{3}{2}$ we have $f_\beta\notin H^2.$ This involves to compute the Taylor series expansion of the function $f(z)=(1-z)^{-1/2}=\sum_{k=0}^\infty c_kz^k.$ By an easy computation we have
\[ c_k=\frac{1}{2^k}\frac{(2k-1)!!}{k!}=\frac{1}{2^k}\frac{(2k)!}{(k!)^2}.\]
Using Stirling's formula we have $c_k\sim \frac{1}{\sqrt{k+1}}.$ Therefore, $f(z)=(1-z)^{-1/2}$ does not belong to $H^2.$\\

Finally we need to deal with the case (iii) $\beta<-\frac{3}{2}.$ For $\beta \leq -2$ the function $f_\beta$ is bounded, and therefore belongs to the Hardy space. We only need to consider the case $-2<\beta <-\frac{3}{2}.$ In such case we write $2+\beta =\frac{1}{2}-\delta,$ where $ 1/2>\delta>0.$
For $0<r<1/2$ the quantity $\|f_\beta(re^{i(\cdot)})\|^2_2$ is uniformly bounded. For $1/2\leq r<1,$ we have
\begin{eqnarray*}
\|f_\beta(re^{i(\cdot)})\|^2_2&=&2\int_0^\pi \frac{1}{\left[{(1-r)^2+4r\sin^2\left(\frac{t}{2}\right)}\right]^{1/2-\delta}}dt\\
&\leq& C \int_0^{\pi/2}\frac{1}{t^{1-2\delta}}dt\\
&\leq& C,
\end{eqnarray*}
where $C$ is a constant independent of $r.$ We thus conclude that $f_\beta\in H^2.$   \qed \par
Below we will give a sequence of functions belonging to different weighted Bergman spaces. Besides their own interest, they will be used to test efficiency of the weight Bergman space POAFD algorithms being applied to functions in different spaces. The reproducing kernel type functions $\frac{1}{(1-\overline{a}z)^{2+\beta}}, a\in \D,$ are all bounded holomorphic functions and hence
in all $\mA_{\alpha},\ \alpha >-1.$

  \section {Experiments}

  \begin{example}
  \end{example}

   Fourier series and POAFD expansion is used to treat the chirp signal $f(t)=\cos t^2.$ The results show that to reach the comparable errors Fourier needs 140 while POAFD needs 32 iterations.

  \begin{figure}[H]
  \centering
  \begin{minipage}[c]{0.5\textwidth}
  \centering
  \includegraphics[height=4.5cm,width=7.5cm]{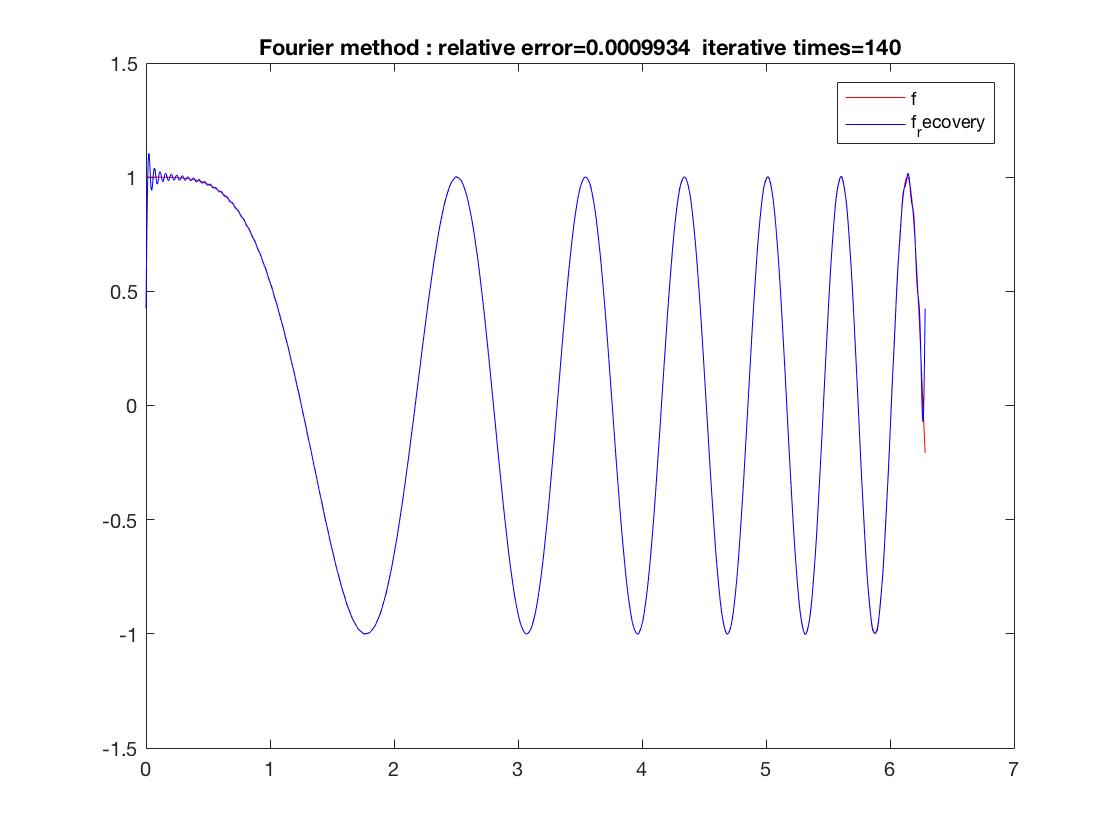}
  \end{minipage}%
  \begin{minipage}[c]{0.5\textwidth}
  \centering
  \includegraphics[height=4.5cm,width=7.5cm]{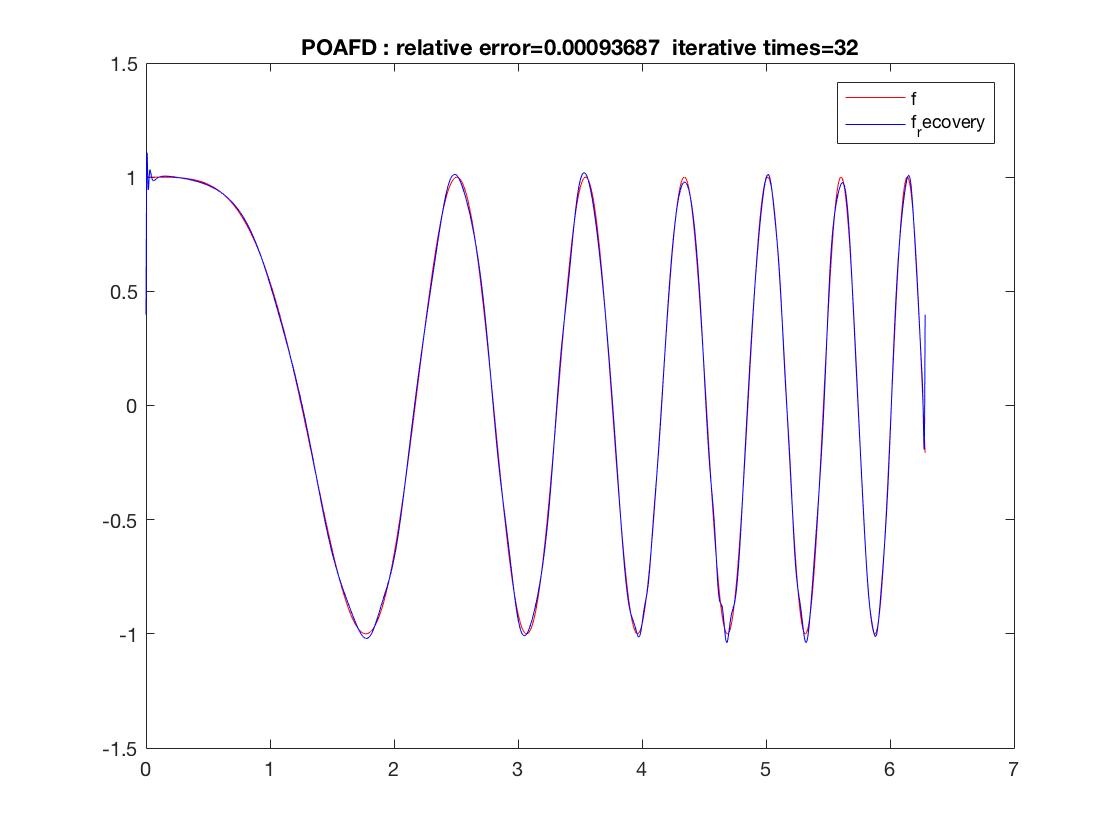}
  \end{minipage}
  \caption*{Figure 9.1}
  \end{figure}
  \begin{figure}[H]

    \centering
    \includegraphics[width=8cm]{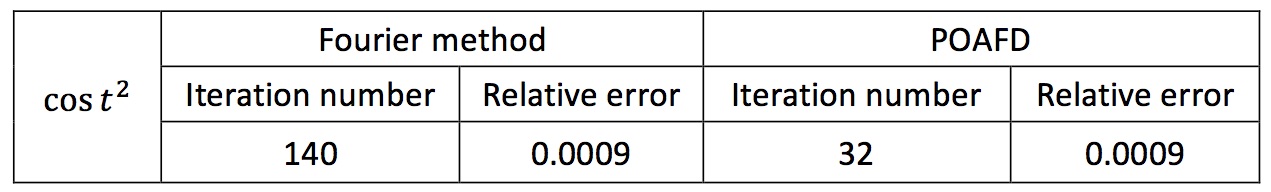}\\
  \caption*{Table 9.1}
  \end{figure}

%
%
%

  \begin{example}
  \end{example}

  The second set of experiments given by figure 9.2 is for recovering functions
  $\sum_{k=0}^{10} \frac{z^k}{(k+1)^{2}}$ and $\prod_{k=1}^{10}\frac{z-a_k}{1-\overline{a}_kz},$ where $a_k$ is randomly selected in $\D.$ The results of 5 iterations are given by figure 9.2. By using the same iteration numbers, POAFD gives smaller relative error than Fourier, $1\times10^{-3}$ and $8\times10^{-6},$ respectively. We also give comparison between different iteration numbers in figure 9.2 with respect to the same relative error. It shows that the iteration number of POAFD is ten iterations less than that of the Fourier method, being more efficient than the latter. The  experiments show that for those examples POAFD outperforms Fourier. \\

  \begin{figure}[H]
  \centering
  \begin{minipage}[c]{0.5\textwidth}
  \centering
  \includegraphics[height=4.5cm,width=7.5cm]{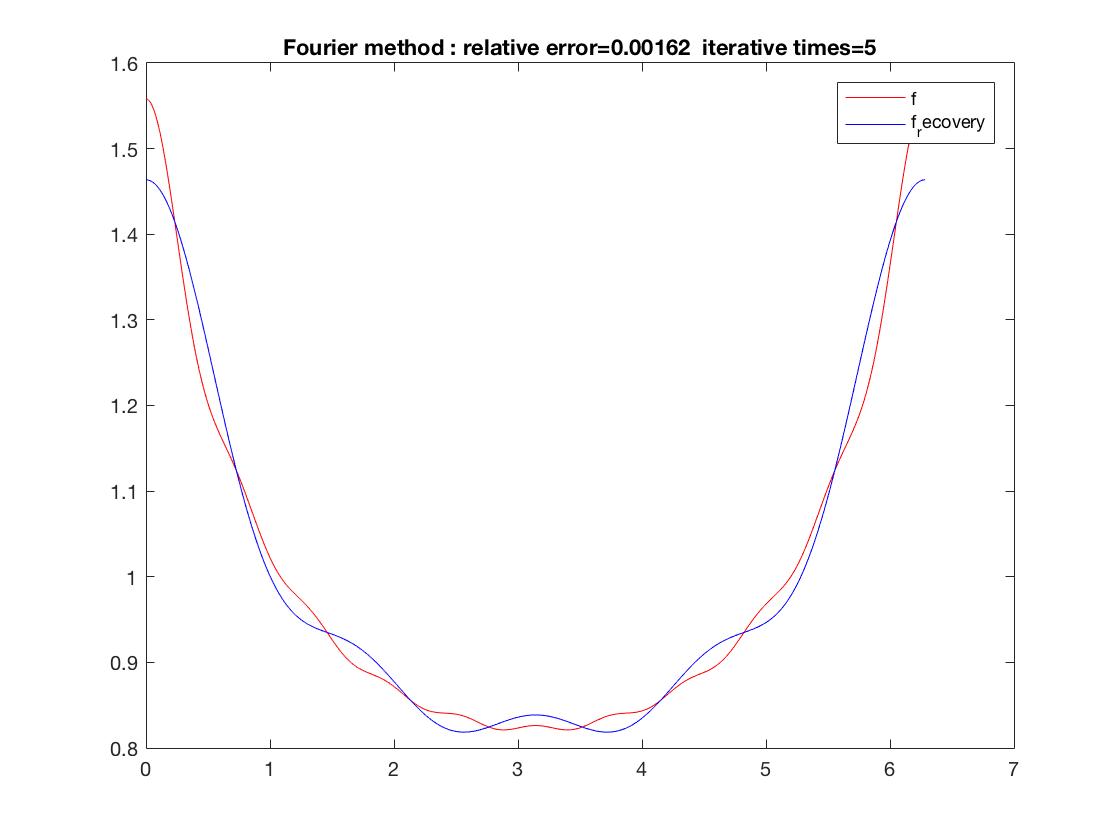}
  \end{minipage}%
  \begin{minipage}[c]{0.5\textwidth}
  \centering
  \includegraphics[height=4.5cm,width=7.5cm]{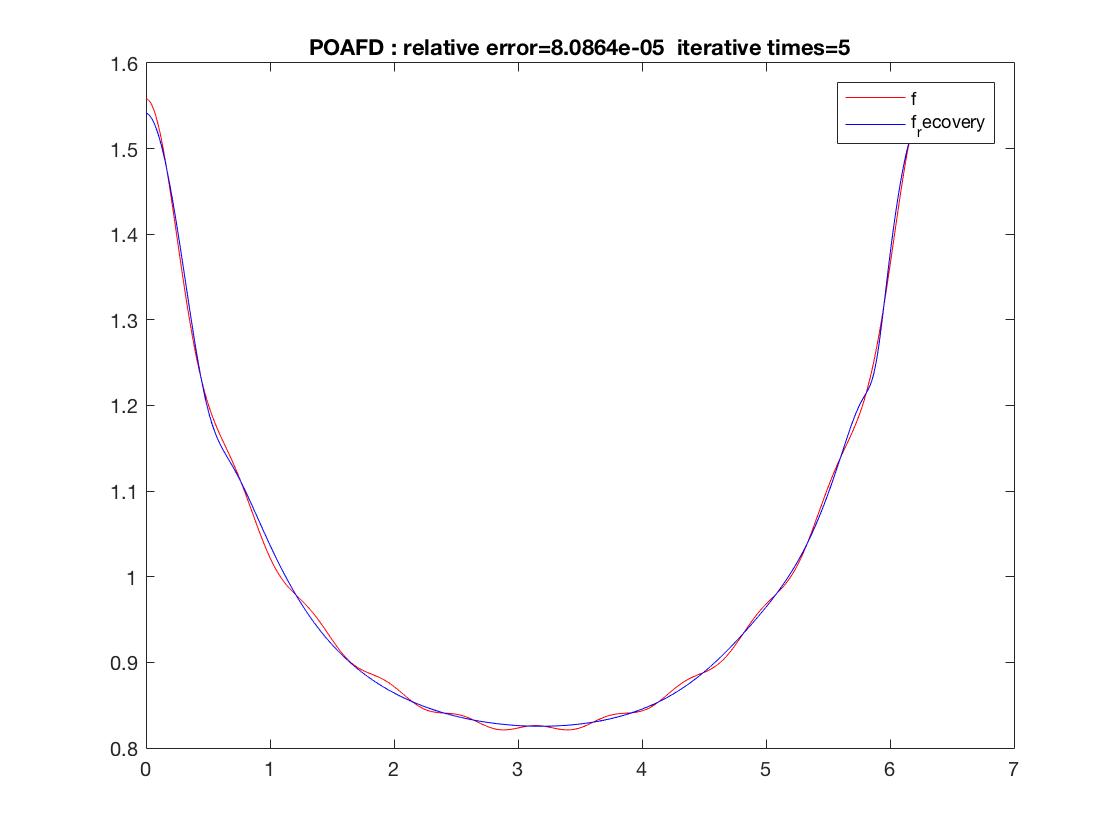}
  \end{minipage}

  \centering
  \begin{minipage}[c]{0.5\textwidth}
  \centering
  \includegraphics[height=4.5cm,width=7.5cm]{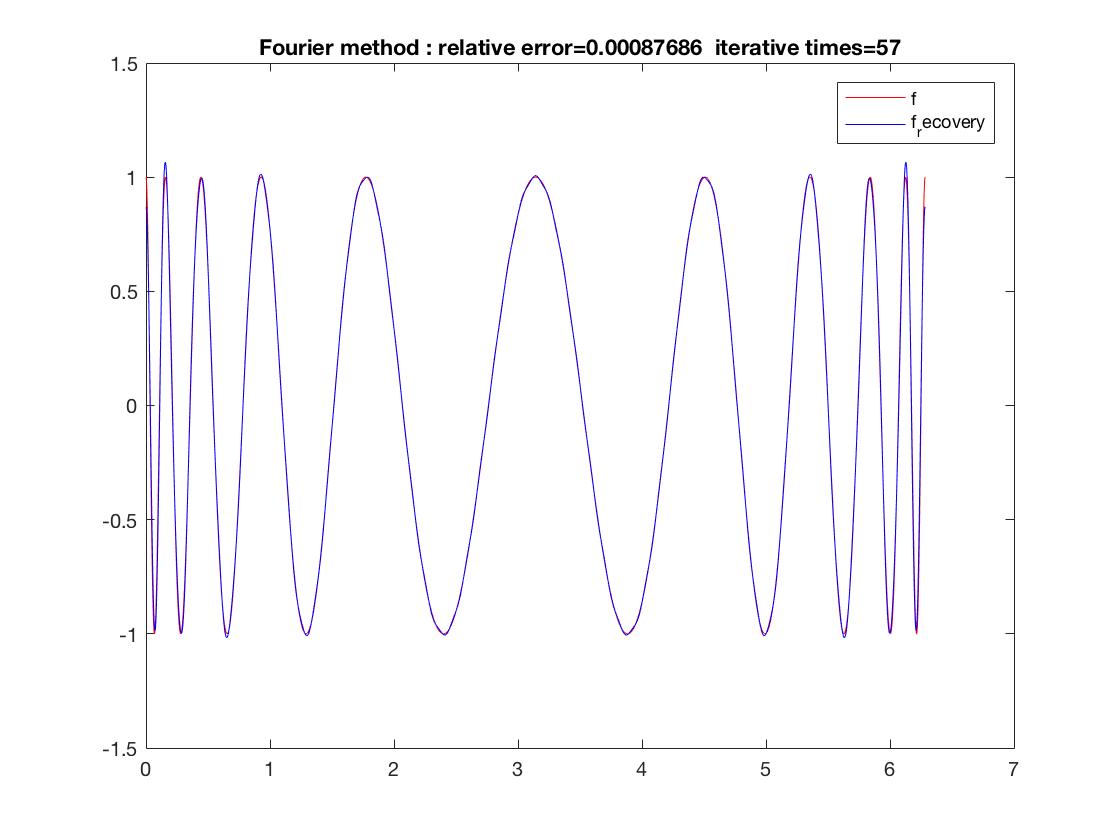}
  \end{minipage}%
  \begin{minipage}[c]{0.5\textwidth}
  \centering
  \includegraphics[height=4.5cm,width=7.5cm]{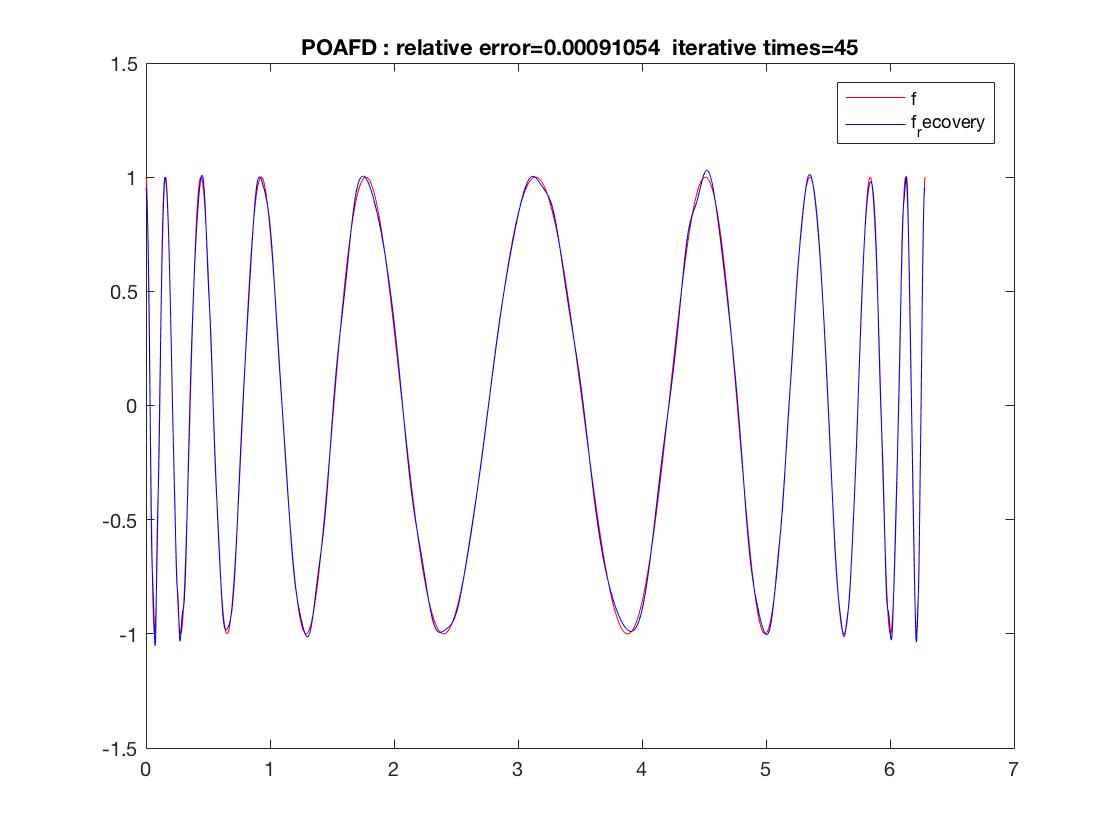}
  \end{minipage}

  \caption*{Figure 9.2}
  \end{figure}

  \begin{figure}[H]

    \centering
    \includegraphics[width=8cm]{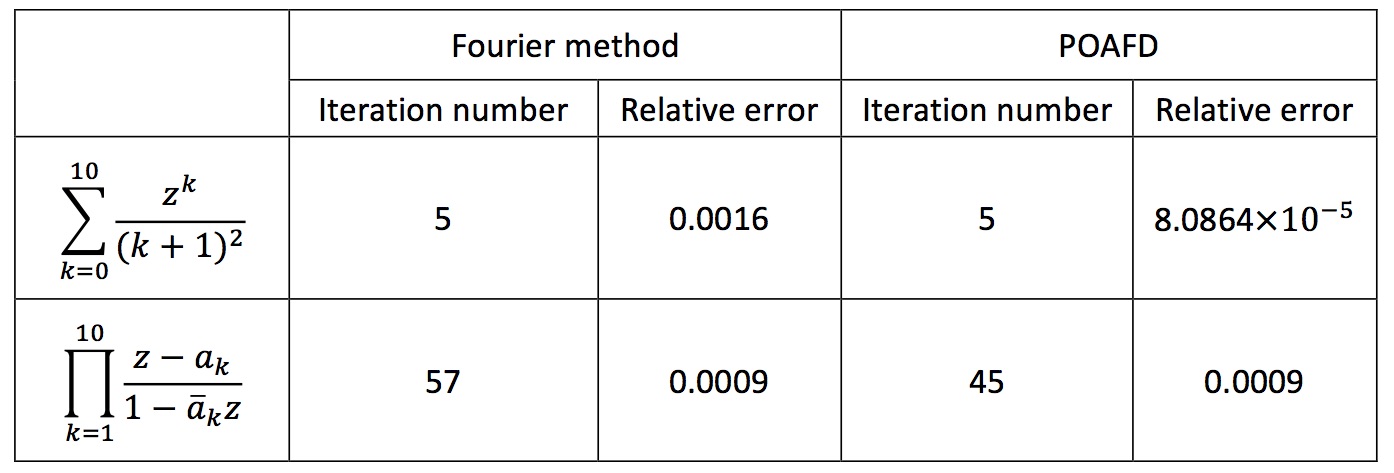}\\
  \caption*{Table 9.2}
  \end{figure}

  \begin{example}
  \end{example}

   The third set of experiments given by figure 9.3 is for recovering a signal with a singular inner function part
  $$\frac{1+2z^2}{(z-2)(z-3)}-\frac{1}{(z+2)(z+3)}e^{z+1+\frac{z+i}{z-i}+\frac{z-i}{z+i}}.$$ Under the same iteration number 150, we have higher convergence speed with smaller error than the Fourier method. It is seen in this example that by POAFD the high oscillatory parts are accurately approximated.\\

  \begin{figure}[H]
  \centering
  \begin{minipage}[c]{0.5\textwidth}
  \centering
  \includegraphics[height=4.5cm,width=7.5cm]{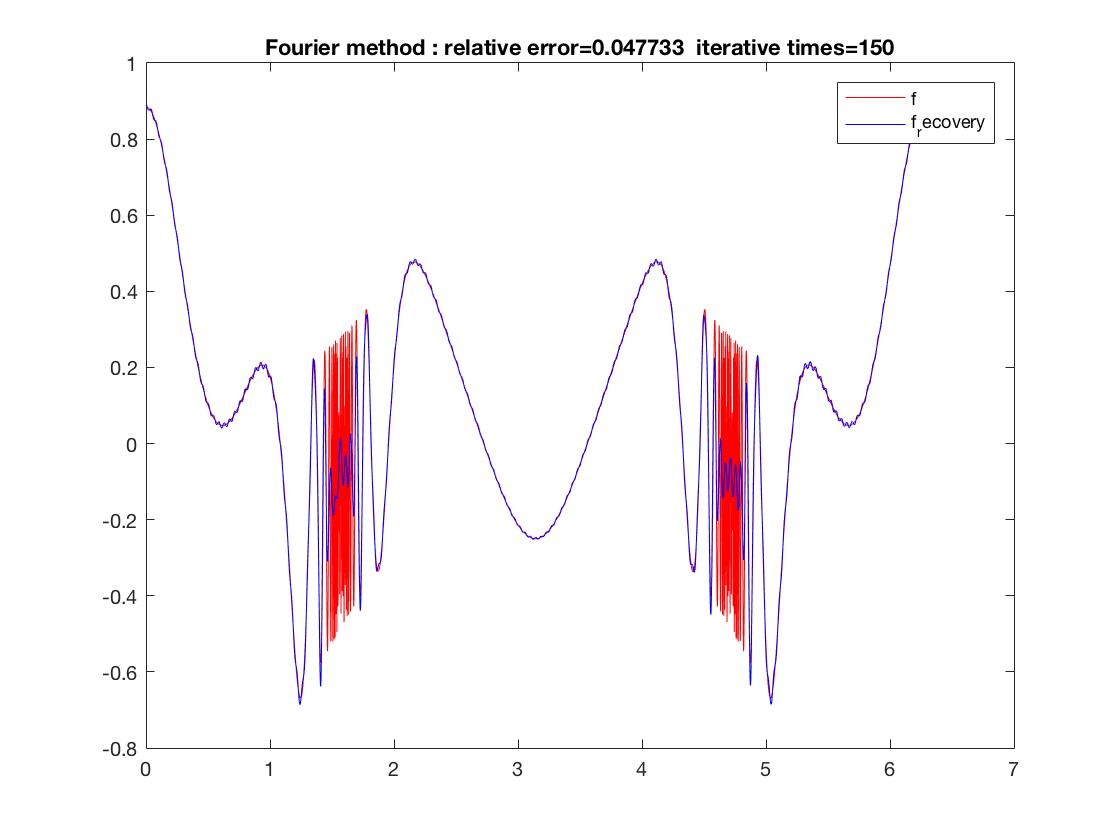}
  \end{minipage}%
  \begin{minipage}[c]{0.5\textwidth}
  \centering
  \includegraphics[height=4.5cm,width=7.5cm]{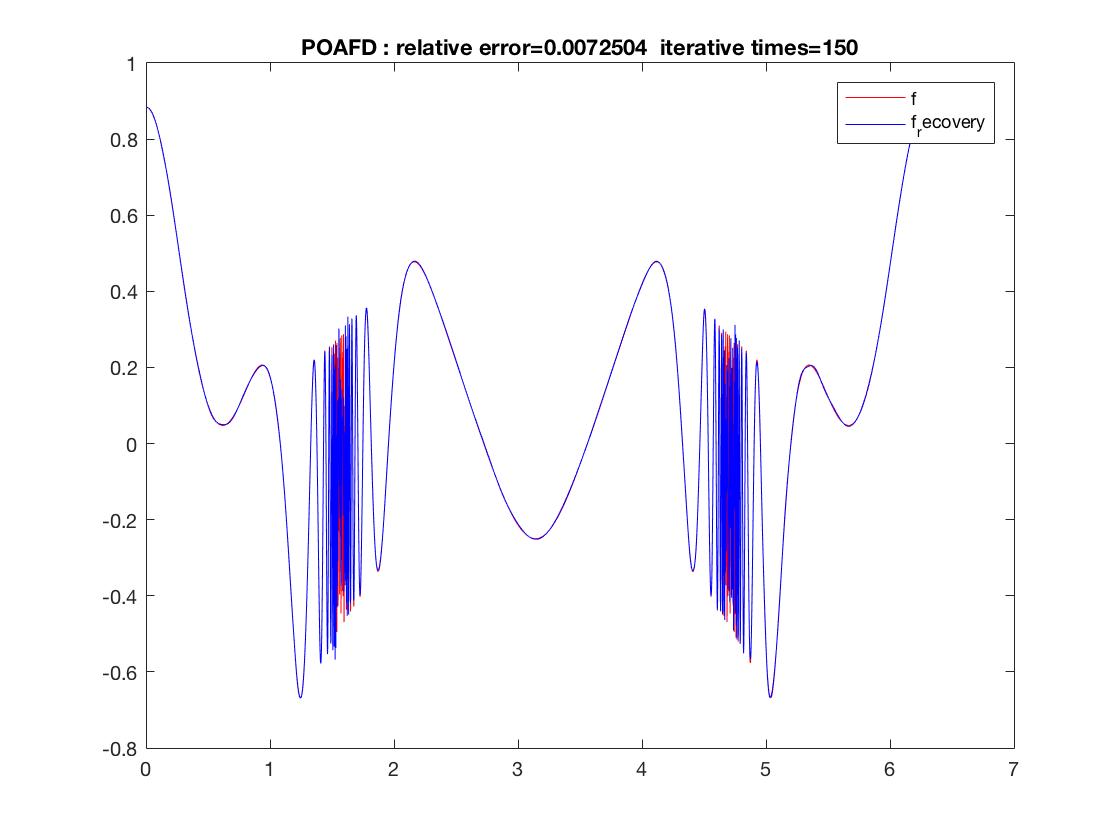}
  \end{minipage}

  \caption*{Figure 9.3}
  \end{figure}\par

  \begin{figure}[H]
    \centering
    \includegraphics[width=12cm]{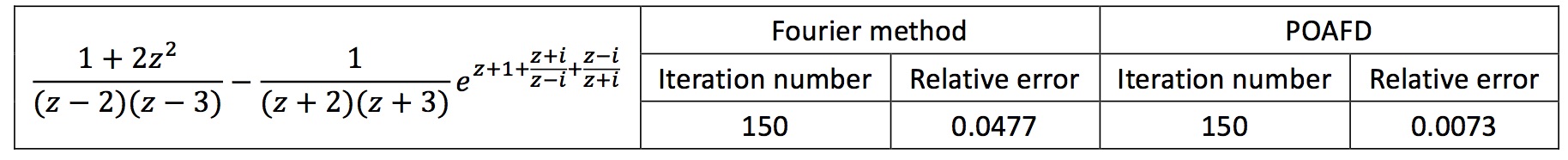}\\
  \caption*{Table 9.3}
  \end{figure}\par

  \begin{example}
  \end{example}

 In this experiment we can summarize from below the enclosed table that for functions of the reproducing kernel type, $\frac{1}{(1-\overline{a}z)^{2+\beta}}, a$ is close to $ \partial\D,$ the POAFD in $\mA_\beta$ performs with the best effect; and when $\alpha$ leaves away from $\beta$  the effectiveness of the POAFD in $\mA_\beta$ for approximating $\frac{1}{(1-\overline{a}z)^{2+\beta}}$ is reduced (IN stands for iteration number, RE for relative error).

  \begin{figure}[H]
    \centering
    \includegraphics[width=12cm]{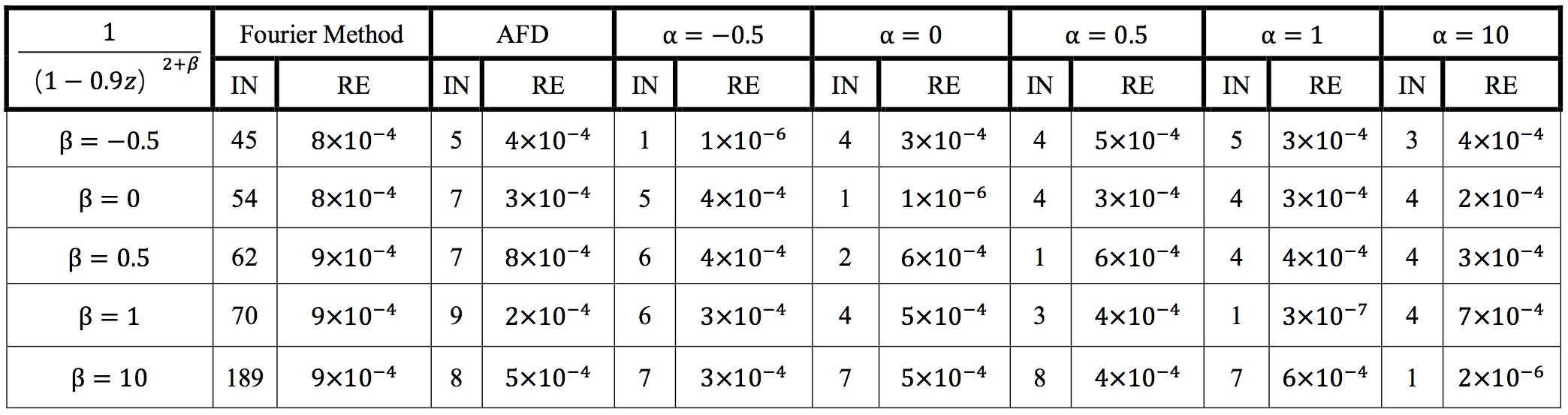}\\
  \caption*{Table 9.4}
  \end{figure}\par

\section{Conclusion}
 The validity of POAFD for the Bergman and weighted Bergman spaces in the unit disc and the upper-half complex plane contexts are evidenced through proving the right boundary vanishing properties or, alternatively, maximal selection theorems. The corresponding algorithms are established. We show that to gain the best approximation at every iterative step occasions of multiple selections of parameters are unavoidable. The process of repeating selection of parameters necessarily induces directional derivatives along with the G-S orthonormalization process. The orthonormal rational systems or weighted Bergman-TM-systems are constructed. The underlying study, in particular, promotes the topic rational approximation to holomorphic functions without boundary limits, and, in particular, with a variety of singularities at the boundary.  The effectiveness and efficiency of the algorithms are tested through different types of Bergman space functions.

\end{document}